%% file: paper.tex
\documentclass[11pt]{article}
\usepackage[top=1in, bottom=1in, left=1in, right=1in]{geometry}
\usepackage{tikz}
\usetikzlibrary{shapes,arrows,matrix,patterns,positioning}
\usepackage{color}
\usepackage{graphicx}
\usepackage{algorithmic}
\usepackage{amssymb}
\usepackage{epstopdf}
\usepackage[norelsize,boxed,linesnumbered,vlined,algo2e]{algorithm2e} 
\usepackage{float}
\usepackage{amsmath,amsthm,bm,color,epsfig,enumerate,caption}
\usepackage{hyperref}
\usepackage{booktabs}
\usepackage{multirow}
\usepackage{array}
\usepackage{bm}

\newtheorem{theorem}{Theorem}[section]

\newtheorem{algo}[theorem]{Algorithm}

\renewcommand{\appendix}[1]{
\section*{Appendix: #1}
}

\renewcommand{\O}{\mathcal{O}}
\renewcommand{\O}{O}

\newcommand{\bbC}{\mathbb{C}}
\newcommand{\bbR}{\mathbb{R}}

\newcommand{\eps}{\epsilon}

\newcommand{\diag}{\text{ diag}}
\newcommand{\coef}[2]{\lambda^{#1}_{#2}}

\makeatletter
\newcommand*{\extendadd}{
  \mathbin{
    \mathpalette\extend@add{}
  }
}
\newcommand*{\extend@add}[2]{
  \ooalign{
    $\m@th#1\leftrightarrow$%
    \vphantom{$\m@th#1\updownarrow$}
    \cr
    \hfil$\m@th#1\updownarrow$\hfil
  }
}
\makeatother

\begin{document}

\title{A Unified Framework for Oscillatory Integral Transform{s}: \\When to use NUFFT or Butterfly Factorization?}

\author{Haizhao Yang
  \vspace{0.1in}\\
  Department of Mathematics, National University of Singapore
}

\maketitle

\begin{abstract}
This paper concerns the fast evaluation of the matvec $g=Kf$ for $K\in \mathbb{C}^{N\times N}$, which is the discretization of an oscillatory integral transform $g(x) = \int K(x,\xi) f(\xi)d\xi$ with a kernel function $K(x,\xi)=\alpha(x,\xi)e^{2\pi\i \Phi(x,\xi)}$, where $\alpha(x,\xi)$ is a smooth amplitude function , and $\Phi(x,\xi)$ is a piecewise smooth phase function with $O(1)$ discontinuous points in $x$ and $\xi$. A unified framework is proposed to compute $Kf$ with $O(N\log N)$ time and memory complexity via the non-uniform fast Fourier transform (NUFFT) or the butterfly factorization (BF), together with an $O(N)$ fast algorithm to determine whether NUFFT or BF is more suitable. This framework works for two cases: 1) explicit formulas for the amplitude and phase functions are known; 2) only indirect access of the amplitude and phase functions are available. Especially in the case of indirect access, our main contributions are: 1) an $O(N\log N)$ algorithm for recovering the amplitude and phase functions is proposed based on a new low-rank matrix recovery algorithm; 2) a new stable and nearly optimal BF with amplitude and phase functions in a form of a low-rank factorization (IBF-MAT) is proposed to evaluate the matvec $Kf$. Numerical results are provided to demonstrate the effectiveness of the proposed framework.
\end{abstract}



\section{Introduction} \label{sec:intro} 
Oscillatory integral {transforms have} been an important topic for scientific computing. After discretization with $N$ grid points in each variable, the integral transform is reduced to a dense matrix-vector multiplication (matvec) $g=Kf$. The direct computation of the matvec takes $O(N^2)$ operations and is prohibitive in large-scale computation. There has been an active research line in developing nearly linear matvec based on the similarity of $K$ to the Fourier matrix \cite{Gang,Alex2}, i.e., $K(x,\xi)=\alpha(x,\xi)e^{2\pi\i  p(x)q(\xi)}$, or based on the complementary low-rank structure of $K$ \cite{Yingwave,Hu,IBF,BF,MBF,Butterfly1,Butterfly2,invRadon,SHT} when the phase function is not in a form of separation of variables. 

The main ideas of existing algorithms are as follows. After computing the low-rank approximation of $\alpha(x,\xi)\approx\sum_{k=1}^r a_k(x)b_k(\xi)$, we have
\begin{equation}
\label{eqn:tg}
g(x)\approx \sum_{k=1}^r a_k(x)\int e^{2\pi\i \Phi(x,\xi)}\left(b_k(\xi)f(\xi)\right)d\xi.
\end{equation}
If $K(x,\xi)=\alpha(x,\xi)e^{2\pi\i  p(x)q(\xi)}$, then 
\[
g(x)\approx \sum_{k=1}^r a_k(x)\int e^{2\pi\i  p(x)q(\xi)}\left(b_k(\xi)f(\xi)\right)d\xi
\]
can be evaluated through $r$ NUFFT's. If the phase function $\Phi(x,\xi)$ is not of the form $p(x)q(\xi)$, then the butterfly factorization (BF) \cite{BF,Butterfly1,Butterfly2} of $e^{2\pi\i \Phi(x,\xi)}$ is computed. The main cost for evaluating \eqref{eqn:tg} is to apply the BF to $r$ vectors, which is $O(rN\log N)$. Hence, after precomputation (low-rank factorization and BF\footnote{In most applications, $K$ is applied to multiple vectors $f$'s. Hence, it is preferable to save the results of expensive computational routines that are independent of the input vectors $f$'s for later applications.}, if needed), both kinds of algorithms admit $O(N\log N)$ computational complexity for applying $K$ to a vector $f$. However, existing algorithms are efficient only when the explicit formulas of the kernel is known (see Table \ref{tab:ex1} and \ref{tab:ex2} for a detailed summary). The computational challenge in the case of indirect access of the kernel function {(see Table \ref{tab:sc} for a list of different scenarios)} motivates a series {of} new algorithms in this paper.

\begin{table}[]
\centering
\begin{tabular}{rccccccc}
\toprule
   Kernels $K(x,\xi)$ & Algorithms & Precomputation time & Application time
                  & memory\\
\toprule
    $\alpha(x,\xi)e^{2\pi\i  p(x)q(\xi)}$ & NUFFT \cite{Gang,Alex2} &   $O(N)$ & $O(N\log N)$ & $O(N)$  \\
    $\alpha(x,\xi)e^{2\pi\i  \Phi(x,\xi)}$ & {\bf NUFFT}  &   $O(N)$ & $O(N\log N)$ & $O(N)$  \\
    $\alpha(x,\xi)e^{2\pi\i  \Phi(x,\xi)}$ & BF \cite{FIO09,IBF} &   $O(N\log N)$ & $O(N\log N)$ & $O(N\log N)$  \\

\bottomrule
\end{tabular}
\caption{Summary of existing algorithms and proposed algorithms (in bold) for the evaluation of $Kf$ when amplitude and phase have explicit formulas. Although the BF in \cite{FIO09} requires no precomputation and $O(N)$ memory, it is a few times slower than the BF in \cite{IBF} regarding the application time. Hence, we adopt the scaling of \cite{IBF} in this paper.}
\label{tab:ex1}
\end{table}

\begin{table}[]
\centering
\begin{tabular}{rccccccc}
\toprule
   Scenarios & Algorithms & Precomputation time & Application time
                  & memory\\
\toprule
      Scenario $1$ & BF \cite{BF} & $O(N^{1.5})$ & $O(N\log N)$ & $O(N\log N)$  \\
 Scenario $2$ & BF \cite{BF} & $O(N^{1.5}\log N)$ & $O(N\log N)$ & $O(N^{1.5})$  \\
 Scenario $3$ & BF \cite{FIO09,IBF} &   $O(N\log N)$ & $O(N\log N)$ & $O(N\log N)$  \\
 All scenarios & {\bf NUFFT} or {\bf IBF-MAT} &    $O(N\log N)$ & $O(N\log N)$ & $O(N\log N)$  \\

\bottomrule
\end{tabular}
\caption{Summary of existing algorithms and proposed algorithms (in bold) for the evaluation of $Kf$ for a general kernel $\alpha(x,\xi)e^{2\pi\i  \Phi(x,\xi)}$ when only indirect access of amplitude and phase is available according to different scenarios {listed in Table \ref{tab:sc}}.}
\label{tab:ex2}
\end{table}

\begin{figure}[ht!]
  \begin{center}
    \begin{tabular}{c}
      \includegraphics[height=2.5in]{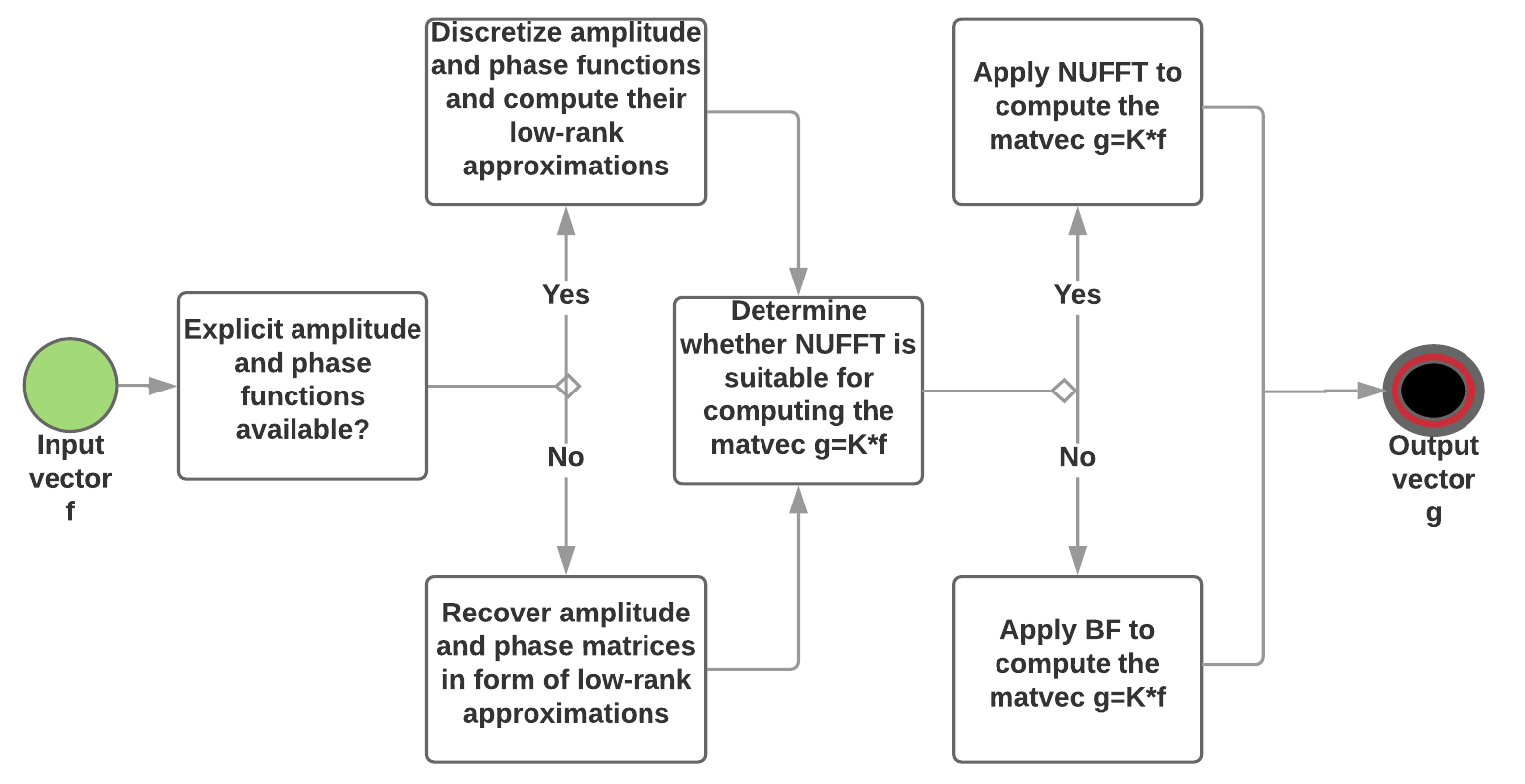}
    \end{tabular}
  \end{center}
  \caption{The computational flowchart of the unified framework using NUFFT or BF. The framework consists of three main steps: 1) construct the low-rank approximations of the amplitude and phase matrices; 2) determine whether NUFFT is applicable; 3) apply NUFFT or BF. When the numerical rank of the phase function $r_\epsilon$ is only larger than the dimension of the problem by one or two, NUFFT is usually faster than BF and hence it will be applied to compute $Kf$. }
\label{fig:flow}
\end{figure}

This paper proposes an $O(N\log N)$ unified framework for evaluating $Kf$ either based on NUFFT or BF (see Figure \ref{fig:flow} for the main computational flowchart of the unified framework). {This framework considers possibly most application scenarios of oscillatory integral transforms. We also briefly discuss how to choose NUFFT or BF to maximize the computational efficiency according to several factors (e.g., accuracy and rank parameters in low-rank factorization, the number of vectors in the matvec) in a serial computational environment.} The unified framework works in two cases: 1) explicit formulas for the amplitude and phase functions are known; 2) only indirect access of the amplitude and phase functions are available. When the explicit formulas are given, computing $Kf$ is relatively simple. Hence, we only focus on the case of indirect access. To the best of our knowledge, the most common indirect access can be summarized into three scenarios in Table \ref{tab:sc}.

\begin{table}[h]
\begin{center}
{\renewcommand\arraystretch{1.25}
\begin{tabular}{|l|l|l|} \hline
Scenario $1$ :& \multicolumn{2}{p{12cm}|}{\raggedright There exists an algorithm for evaluating an arbitrary entry of the kernel matrix in $O(1)$ operations \cite{James:2017,Bremer201815,BF,Butterfly2}.} \\ \hline
Scenario $2$ :& \multicolumn{2}{p{12cm}|}{\raggedright There exist an $O(N\log N)$ algorithm for applying $K$ and its transpose to a vector \cite{LUBF,BF,HSSBF,precon2}.} \\ \hline
Scenario $3$ :& \multicolumn{2}{p{12cm}|}{\raggedright The amplitude and the phase functions are solutions of partial differential equations (PDE's) \cite{Yingwave}. $O(1)$ columns and rows of the amplitude and phase matrices are available by solving PDE's.} \\ \hline
\end{tabular}}
\end{center}
\caption{Three scenarios of the indirect access of the amplitude and phase functions. 
}
\label{tab:sc}
\end{table}

As the first main contribution of this paper, in the case of indirect access, a nearly linear scaling algorithm is proposed to recover the amplitude and phase matrices in a form of low-rank matrix factorization. As far as we know, the low-rank matrix recovery problem in this paper has not been studied before since there is no direct access to the entries of low-rank matrices. Hence, there is no existing algorithm in the literature suitable for this problem.

As the second main contribution, when the low-rank amplitude and phase matrices have been recovered, a new BF (named as IBF-MAT for short) is proposed for the matvec $Kf$. IBF-MAT is the first BF for the matvec $Kf$ with $O(N\log N)$ complexity for both precomputation and application in the case of indirect access (see Table \ref{tab:ex2} for the comparison with existing algorithms).

Finally, this paper shows that if the numerical rank of $\Phi(x,\xi)$ is $r_\epsilon$ (depending on an $\epsilon$ accuracy parameter), a $r_\epsilon$-dimensional NUFFT can be applied to evaluate \eqref{eqn:tg} in $O(N\log N)$ operations. The dimension of the NUFFT, $r_\epsilon$, could be larger than the dimension of the variables $x$ and $\xi$, and hence we consider it as a dimension lifting technique. This new method significantly extend the application range of the NUFFT approach for computing $Kf$.

The rest of the paper is organized as follows. In Section \ref{sec:pre}, we revisit existing low-rank factorization techniques and propose our new low-rank matrix factorization in the case of indirect access. In Section \ref{sec:frame}, we introduce the new NUFFT approach by dimension lifting. In Section \ref{sec:BF},  we introduce the IBF-MAT. Finally, we provide several numerical examples to demonstrate the efficiency of the proposed unified framework in Section \ref{sec:results}.

%======================================================
\section{Low-rank matrix factorization}
\label{sec:pre}
This section is for the first main step in the unified framework as shown in Figure \ref{fig:flow}: low-rank matrix factorizations of the amplitude and phase matrices.

\subsection{Existing low-rank matrix factorization}
\label{sec:LRF}

%Two classical techniques for low-rank matrix factorizations are introduced here. 

%\vspace{0.25cm}
\textbf{Low-rank approximation by randomized sampling}

\vspace{0.25cm}
For $K \in \bbC^{m\times n}$, we define a rank-$r$
approximate singular value decomposition (SVD) of $K$ as
\begin{equation}\label{eqn:SVD}
K \approx U_0 \Sigma_0 V_0^*,
\end{equation}
where $U_0\in \bbC^{m\times r}$ is orthogonal, $\Sigma_0\in \bbR^{r\times
  r}$ is diagonal, and $V_0\in \bbC^{n\times r}$ is orthogonal. Efficient randomized tools have been proposed to compute approximate SVDs for numerically low-rank matrices \cite{Butterfly4,randomSVD}. The one  in \cite{Butterfly4} (see Algorithm \ref{alg:rlr}) is more attractive because it only requires $\O(m+n)$ operations and memory. We adopt MATLAB notations to describe Algorithm \ref{alg:rlr} for simplicity: given row and column
index sets $I$ and $J$, $K_{I,J}=K(I,J)$ is the
submatrix with entries from rows in $I$ and columns in $J$; the index set for an entire row or column is denoted as $``:"$.

\vspace{0.25cm}
\begin{algorithm2e}[H]
\caption{Randomized sampling for a rank-$r$ approximate SVD.}
\label{alg:rlr}

 Input: A matrix $K \in \bbC^{m\times n}$, a rank parameter $r$, and an over-sampling parameter $q$.

Output: The low-rank factorization $K\approx U_0\Sigma_0V_0^*$ in \eqref{eqn:SVD}.

  Let $\Pi_{col}$ and $\Pi_{row}$ denote the important columns and rows of $K$
  that are used to form the column and row bases.
  Initially $\Pi_{col} = \emptyset$ and $\Pi_{row} = \emptyset$.

  Randomly sample $rq$ rows and denote their indices by $S_{row}$.  Let $I
  = S_{row}\cup \Pi_{row}$. Perform a pivoted {QR} decomposition of $K_{I,:}$ to get  $K_{I,:} P = QR$,  where $P$ is the resulting permutation matrix and $R=(r_{ij})$ is an
  $\O(r)\times n$ upper triangular matrix. Define the important column
  index set $\Pi_{col}$ to be the first $r$ columns picked within the
  pivoted QR decomposition.

  Randomly sample $rq$ columns and denote their indices by $S_{col}$.
  Let $J = S_{col}\cup \Pi_{col}$. Perform a pivoted {LQ} decomposition of
  $K_{:,J}$ to get $P K_{:,J} = LQ$, where $P$ is the resulting permutation matrix and $L=(l_{ij})$ is an
  $m\times \O(r)$ lower triangular matrix.  Define the important
  row index set $\Pi_{row}$ to be the first $r$ rows picked within the
  pivoted LQ decomposition.

  Repeat Line $4$ and $5$ a few times to ensure $\Pi_{col}$ and $\Pi_{row}$
  sufficiently sample the important columns and rows of $K$.

  Apply the pivoted QR factorization to $K_{:,\Pi_{col}}$ and let
  $Q_{col}$ be the matrix of the first $r$ columns of the $Q$ matrix.
  Similarly, apply the pivoted QR factorization to $K_{\Pi_{row},:}^*$
  and let $Q_{row}$ be the matrix of the first $r$ columns of the $Q$
  matrix.

  Let $S_{col}$ and $S_{row}$ be the index sets of a few extra randomly
  sampled columns and rows.  Let $J = \Pi_{col} \cup S_{col}$ and $I =
  \Pi_{row} \cup S_{row}$.  Let $M = (Q_{col})_{I,:}^\dagger K_{I,J}
  (Q_{row}^*)_{:,J}^\dagger$, where $(\cdot)^\dagger$ stands for the
  pseudo-inverse.

  Compute the {SVD} $M \approx U_M\Sigma_MV_M^*$ and let $U_0 = Q_{col} U_M$, $\Sigma_0 = \Sigma_M$, and $V_0^* = V_M^*Q_{row}^* $.
\end{algorithm2e}
\vspace{1cm}

\textbf{Interpolative low-rank approximation}

\vspace{0.25cm}
Algorithm \ref{alg:rlr} is sufficiently efficient if we allow a linear complexity to construct the low-rank approximation. However, to construct the BF in nearly linear operations, we cannot even afford linear scaling low-rank approximations; we can only afford an algorithm that provides the low-rank factors with explicit formulas. This motivates the interpolative low-rank approximation below.

Let us focus on the case of a kernel function $K(x,\xi)=e^{2\pi\i  \Phi(x,\xi)}$ and its discretization $K=e^{2\pi\i \Phi}\in \mathbb{C}^{N_A\times N_B}$ to introduce the interpolative low-rank approximation. We assume that $x$ and $\xi$ are one-dimensional variables and the algorithm below can be easily generalized to higher dimensional cases by tensor products. Note that if the phase function is given in a form of separation of variables, i.e., $\Phi(x,\xi)=\sum_{k=1}^r u_k(x)v_k(\xi)$, the following interpolative factorization will also work with a minor modification. 

Let $A$ and $B$ denote the sets of {contiguous} row and column indices of $K$. If $A\times B$ corresponds to a small two-dimensional interval in the variables $x\times \xi$, then a low-rank approximation
\[
K(A,B)=e^{2\pi \i \Phi(A,B)} \approx U_0 V_0^*
\]
exists and can be constructed via Lagrange interpolation as follows.

Suppose the numbers of elements in $A$ and $B$ are $N_A$ and $N_B$, respectively.
Let
\begin{equation}
  R(A,B) := \Phi(A,B)-ones(N_A,1)*\Phi(c_A,B)-\Phi(A,c_B)*ones(1,N_B) +
  \Phi(c_A,c_B),
  \label{eqn:RAB}
\end{equation}
where {$c_A$ and $c_B$ are the indices of $A$ and $B$ closest to the mean of all indices in $A$ and $B$}, 
respectively, then $K$ can be written as
\begin{equation}
  K(A,B)=e^{-2\pi\i\Phi(c_A,c_B)}*\diag\left(e^{2\pi\i \Phi(A,c_B)}\right)*
 e^{2\pi\i R(A,B)}  *\diag\left(e^{2\pi\i \Phi(c_A,B)}\right).
  \label{eqn:terms}
\end{equation}
Hence, {the low-rank approximation of $e^{2\pi\i R(A,B)}$ immediately gives
the low-rank approximation of $K(A,B)$}.
A Lagrange interpolation can be applied to construct the
low-rank approximation of $e^{2\pi\i R(A,B)}$. 

\newcommand{\talpha}{\tilde{\alpha}}
\newcommand{\tbeta}{\tilde{\beta}}

Recall the challenge that we may not have explicit formulas for the amplitude or phase functions. Hence, we cannot use Chebyshev grid points in the Lagrange interpolation to maintain a small uniform error as the previous BF in \cite{FIO09,IBF} does. Therefore, we choose indices in $A$ or $B$ in a similar manner like {Mock-Chebyshev points} \cite{Mock2,Mock1} as follows\footnote{
{Though it was shown in \cite{Impossibility} that no fast stable approximation of analytic functions from equispaced samples in a bounded interval in the sense of $L^\infty$-norm with an exponential convergence rate is available, the Mock-Chebyshev points admit polynomial interpolation with a root-exponential convergence rate. In this paper, we care more about the approximation error at the equispaced sampling locations, in which case it is still unknown whether the Mock-Chebyshev points admit an exponential convergence rate. }}.

Let us assume $A=\{1,\dots,N_A\}$ and $B=\{1,\dots,N_B\}$. If an index set doesn't start with the index $1$, we can simply shift the grid points accordingly. For a fixed integer $r$, the Chebyshev grid of order $r$ on $[-\frac{1}{2},\frac{1}{2}]$ is defined by 
\[
\left\{ z_t=\frac{1}{2}\cos\left(\frac{(t-1)\pi}{r-1}\right)\right\}_{1\leq t\leq r}.
\]
A grid adapted to the index set $A$ is then defined via shifting, scaling, and rounding as
\begin{equation}
\label{eqn:gdA}
\{x_t\}_{t=1,\dots,r}=\left\{\text{Round}\left(t+(N_A-r) (z_t+\frac{1}{2})\right)\right\}_{t=1,\dots,r}.
\end{equation}
Note that the rounding operator may result in repeated grid points. Only one grid point will be kept if repeated. Similarly, a grid adapted to the index set $B$ is defined as
\begin{equation}
\label{eqn:gdB}
\{\xi_t\}_{t=1,\dots,r}=\left\{\text{Round}\left(t+(N_B-r) (z_t+\frac{1}{2})\right)\right\}_{t=1,\dots,r}.
\end{equation}

Given a set of indices $\{x_t\}_{t=1,\ldots, r}$ in $A$, define
Lagrange interpolation polynomials $M^A_t(x)$ by
\[
M^A_t(x)=\prod_{1\leq j\leq r,j\neq t} \frac{x-x_j}{x_t-x_j}.
\]
Similarly, $M_t^B$ is denoted as the Lagrange interpolation
polynomials for $B$.

Now we are ready to construct the low-rank approximation of $e^{2\pi\i R(A,B)}$ by interpolation:
\begin{itemize}
  \item when we interpolate in $\xi$, the low-rank approximation of
    $e^{2\pi\i R(A,B)}$ is given by
    \begin{equation}
   e^{2\pi\i R(A,B)} \approx U_0 V_0^*,
      \label{eqn:intp1}
    \end{equation}
    where \[
   U_0= \begin{pmatrix}
    e^{2\pi\i R(A,\xi_1)},\hdots,e^{2\pi\i R(A,\xi_r)}
    \end{pmatrix}\in \mathbb{C}^{N_A\times r},
    \]
    \[
   V_0= \begin{pmatrix}
     (M^B_1(B))^*,\hdots,(M^B_r(B))^*
     \end{pmatrix}\in \mathbb{C}^{N_B\times r},
    \]
    and each $M^B_t(B)$ denotes a row vector of length $N_B$ such that the $k$-th entry  is
    \[
    M^B_t(\xi_k)=\prod_{1\leq j\leq r,j\neq t} \frac{\xi_k-\xi_j}{\xi_t-\xi_j}
    \]
     for $\xi_k\in B$, $k=1,\dots,N_B$, given by \eqref{eqn:gdB}.
  \item when we interpolate in $x$, the low-rank approximation of
    $e^{2\pi\i R(A,B)}$  is
     \begin{equation}
   e^{2\pi\i R(A,B)} \approx U_0 V_0^*,
      \label{eqn:intp2}
    \end{equation}
    where
    \[
   U_0= \begin{pmatrix}
     (M^A_1(A))^*,\hdots,(M^A_r(A))^*
     \end{pmatrix}\in  \mathbb{C}^{N_A\times r},
    \] 
    \[
   V_0= \begin{pmatrix}
    \left(e^{2\pi\i R(x_1,B)}\right)^*,\hdots,\left(e^{2\pi\i R(x_r,B)}\right)^*
    \end{pmatrix}\in \mathbb{C}^{N_B\times r},
    \]
    and each $M^A_t(A)$ denotes a row vector of length $N_A$ such that the $k$-th entry is
    \[
    M^A_t(x_k)=\prod_{1\leq j\leq r,j\neq t} \frac{x_k-x_j}{x_t-x_j}
    \]
    for $x_k\in A$, $k=1,\dots,N_A$, given by \eqref{eqn:gdA}.
  \end{itemize}

Finally, we are ready to construct the low-rank approximation for the
matrix $e^{2\pi\i \Phi(A,B)}$ when we have $\Phi(A,B)$ or equivalently a low-rank factorization of $\Phi(A,B)$ as in Algorithm \ref{alg:ilr}.

\begin{algorithm2e}[]
\caption{Interpolative low-rank approximation for one-dimensional kernel $e^{2\pi\i \Phi(x,\xi)}$. Factorization in higher dimensions can be constructed similarly via tensor products.}
\label{alg:ilr}
 Input: The phase matrix $\Phi\in \mathbb{C}^{N\times N}$ or its low-rank factorization {$\Phi=\bar{U}\bar{V}^*$}. Contiguous index sets $A$ and $B$ of the row and column indices of $\Phi$, respectively. A rank parameter $r$.

Output: The low-rank factorization $UV^*$ such that $UV^*\approx e^{2\pi\i \Phi(A,B)}$, where $U\in\mathbb{C}^{N_A\times r}$, and $V\in\mathbb{C}^{N_B\times r}$, where $N_A$ is the number of elements in $A$ and $N_B$ is for $B$.

\If{the input contains low-rank factors $\bar{U}$ and $\bar{V}$ of $\Phi$}{define a function to evaluate an arbitrary entry of $\Phi$ at the position $(m,n)$ in $O(1)$ operations as follows
\[
\Phi(m,n) = \bar{U}(m,:)\bar{V}(n,:)^*.
\]}

\If{interpolation in the variable $\xi$ in $B$}{by \eqref{eqn:terms} and \eqref{eqn:intp1}, we have
\begin{equation}
\label{eqn:ab1}
U:= e^{-2\pi\i\Phi(c_A,c_B)}*\diag\left(e^{2\pi\i \Phi(A,c_B)}\right)*
U_0,\quad V^*:=V_0^* *\diag\left(e^{2\pi\i \Phi(c_A,B)}\right),
\end{equation}
where $U_0$ and $V_0$ are given just below \eqref{eqn:intp1}.}
\If{interpolation in the variable $x$ in $A$}{by \eqref{eqn:terms} and \eqref{eqn:intp2}, we have
\begin{equation}
\label{eqn:ab2}
U:= e^{-2\pi\i\Phi(c_A,c_B)}*\diag\left(e^{2\pi\i \Phi(A,c_B)}\right)*
 U_0,\quad V^*:=V_0^* *\diag\left(e^{2\pi\i \Phi(c_A,B)}\right),
\end{equation}
where $U_0$ and $V_0$ are given just below \eqref{eqn:intp2}.}
\end{algorithm2e}

\vspace{0.25cm}

\subsection{New low-rank matrix factorization with indirect access}
\label{sec:nLR}

This section introduces a nearly linear scaling algorithm for constructing the low-rank factorization of the phase matrix $\Phi\in\mathbb{R}^{N\times N}$ when we only know the kernel matrix $K=e^{2\pi\i  \Phi}$ through Scenarios $1$ and $2$ in Table \ref{tab:sc}. The main idea is to recover $O(1)$ randomly selected columns and rows of $\Phi$ from the corresponding columns and rows of $K=e^{2\pi \i\Phi}$. Then by Algorithm \ref{alg:rlr} in Section \ref{sec:LRF}, we can construct the low-rank factorization of $\Phi$.

Obtaining $O(1)$ randomly selected columns and rows of $K$ is simple in Scenarios $1$ and $2$: we can directly evaluate them in Scenario $1$; we apply the kernel matrix $K$ and its transpose to $O(1)$ randomly selected natural basis vectors in $\mathbb{R}^N$ to obtain the {columns} and rows. 

However, reconstructing the corresponding columns and rows of $\Phi$ from those of $K=e^{2\pi\i\Phi}$ is more {challenging}. The difficulty comes from the fact that
\[
\frac{1}{2\pi} \Im\left(\log \left(K(i,j)\right)\right) =\frac{1}{2\pi}\Im\left(\log\left(e^{2\pi\i \Phi(i,j)}\right)\right) =  \frac{1}{2\pi}\arg\left( e^{2\pi\i \Phi(i,j)}\right)= \mod(\Phi(i,j),1),
\]
where $\Im(\cdot)$ returns the imaginary part of the complex number, and $\arg(\cdot)$ returns the argument of a complex number. Hence, $\Phi$ is only known up to modular $1$. 

 Fortunately, our main purpose is not to recover the exact $\Phi$ that generates $K$; instead, we are interested in a low-rank matrix $\Psi$ such that 
 \begin{equation}\label{eqn:tr}
 \mod(\Psi,1) = \frac{1}{2\pi} \Im\left(\log \left(K\right)\right).
 \end{equation} 
  Based on the smoothness of the phase function, a $TV^3$-norm\footnote{{The $TV^3$-norm of a vector $v\in\mathbb{R}^N$ is defined as $\|v\|_{TV^3}:= \sum_{i=2}^{N-2} |v_{i+1}+v_{i-1}-2v_i- (v_{i+2}+v_{i}-2v_{i+1})|$ in this paper. Similarly, The $TV^1$-norm of a vector $v\in\mathbb{R}^N$ is defined as $\|v\|_{TV^1}:= \sum_{i=2}^{N} |v_{i}-v_{i-1}|$. The $TV^2$-norm of a vector $v\in\mathbb{R}^N$ is defined as $\|v\|_{TV^2}:= \sum_{i=2}^{N-1} |v_{i+1}+v_{i-1}-2v_i|$. }} minimization technique is proposed to recover the columns and rows of $\Phi$ up to an additive error matrix $E$ that is numerically low-rank, i.e., the $TV^3$-norm minimization technique returns a matrix $\Psi=\Phi+E$ such that $e^{2\pi\i \Psi}=e^{2\pi\i \Phi}$ and $E$ is numerically low-rank.  
 
  To be more rigorous, we look for the solution of the following combinatorial constrained $TV^3$-norm minimization problem: 
  \begin{eqnarray}
 \label{eqn:mintv}
& \smash{\displaystyle\min_{\Phi\in\mathbb{R}^{N\times N}}} & \sum_{i\in\mathcal{R}}\|\Phi(i,:)\|_{TV^3}+ \sum_{j\in\mathcal{C}}\|\Phi(:,j)\|_{TV^3}\\
& \text{subject to} &  \mod(\Phi(i,j),1) = \frac{1}{2\pi} \Im\left(\log \left(K(i,j)\right)\right)\nonumber\\
& & \text{ for }i\in \mathcal{R} \text{ or } j\in\mathcal{C},\nonumber
 \end{eqnarray}
 where $\mathcal{C}$ and $\mathcal{R}$ are column and row index sets with $O(1)$ randomly selected indices, respectively.
 
 {The problem addressed here is similar to phase retrieval problems, but has a different setting to existing phase retrieval applications and different aims in numerical computation. Phase retrieval problems usually have sparsity assumptions on the signals (or after an appropriate transformation) that lose phases. In the problem considered in this paper, $e^{2\pi\i \Phi}$ is dense and might not be sparse after a transformation (e.g., the Fourier transform or wavelet transform). Furthermore, there are only $O(N)$ samples of the target matrix of size $N\times N$ to be recovered and the hard constrain \eqref{eqn:mintv} is preferred instead of treating it as a soft constrain. $TV^1$-norm is a useful tool for regularization in phase retrieval problems; however, $TV^3$-norm is preferred in this paper since, for example, $\{\Phi(x,y)+ax+by\}_{a,b\in \mathbb{Z}}$ are good solutions to obtain the low-rank factorization of the phase function, and it is not necessary to pick up one function among $\{\Phi(x,y)+ax+by\}_{a,b\in \mathbb{Z}}$ with the minimum $TV^1$-norm using much extra effort. $TV^3$-norm minimization leave us much more flexibility to obtain an approximately good solution to \eqref{eqn:tr} quickly.}
 
 {Our goal here is an $O(N)$ algorithm for solving the matrix recovery problem in \eqref{eqn:mintv}. Though there have been many efficient algorithms for phase retrieval problems, they usually require computational cost at least $O(nN^2)$, where $N^2$ is the size of the target and $n$ is the number of iterations. $n$ and $N^2$ are both too large to be applied in our problem. Hence, instead of solving \eqref{eqn:mintv} exactly using advanced optimization techniques, we propose a heuristic fast algorithm to identify a reasonably good approximate solution to \eqref{eqn:mintv}. As we can see in numerical examples, the proposed heuristic algorithm works well in most applications.} 
 
\begin{figure}[ht!]
  \begin{center}
    \begin{tabular}{ccccc}
      \includegraphics[height=1in]{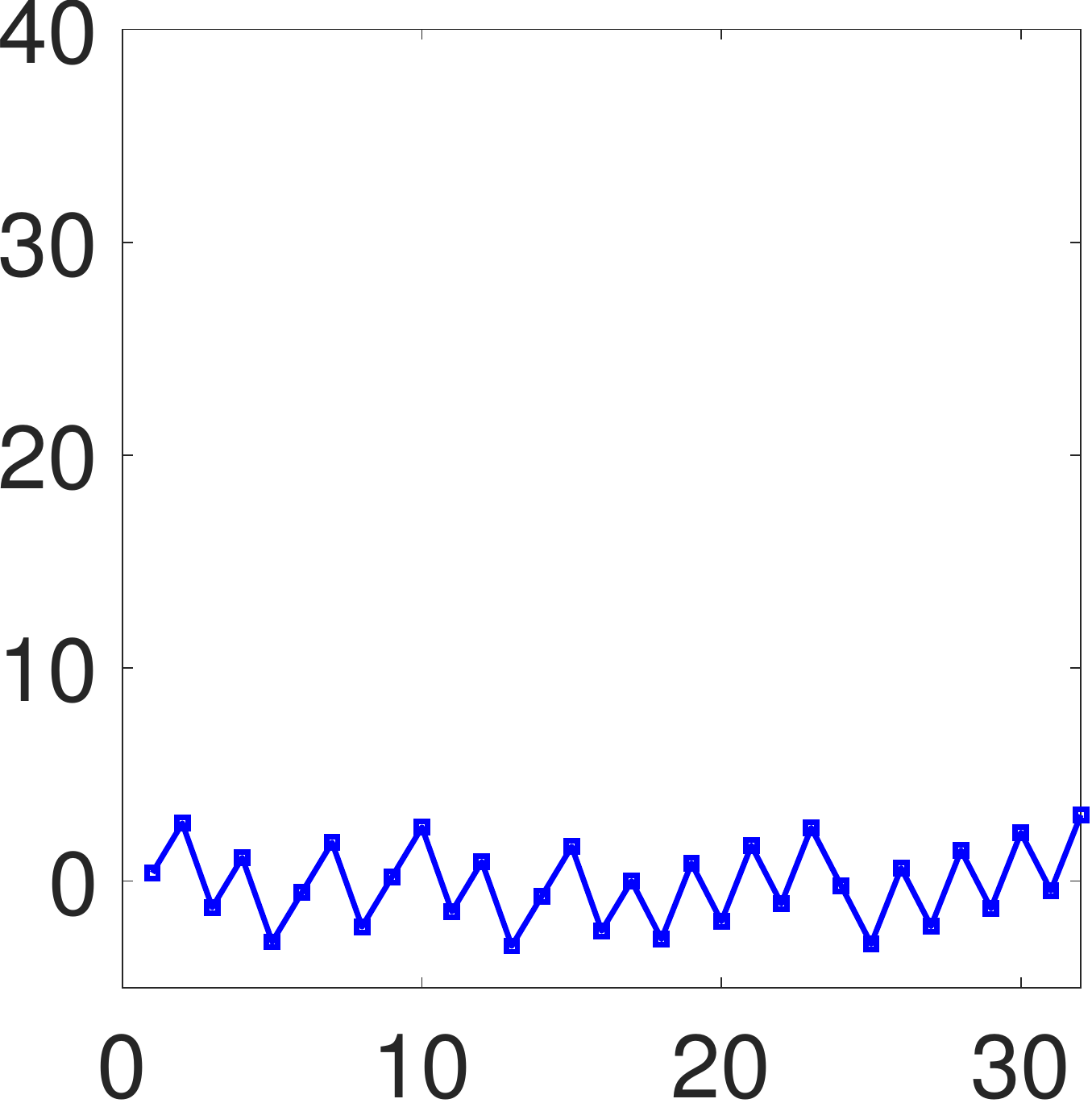}&\includegraphics[height=1in]{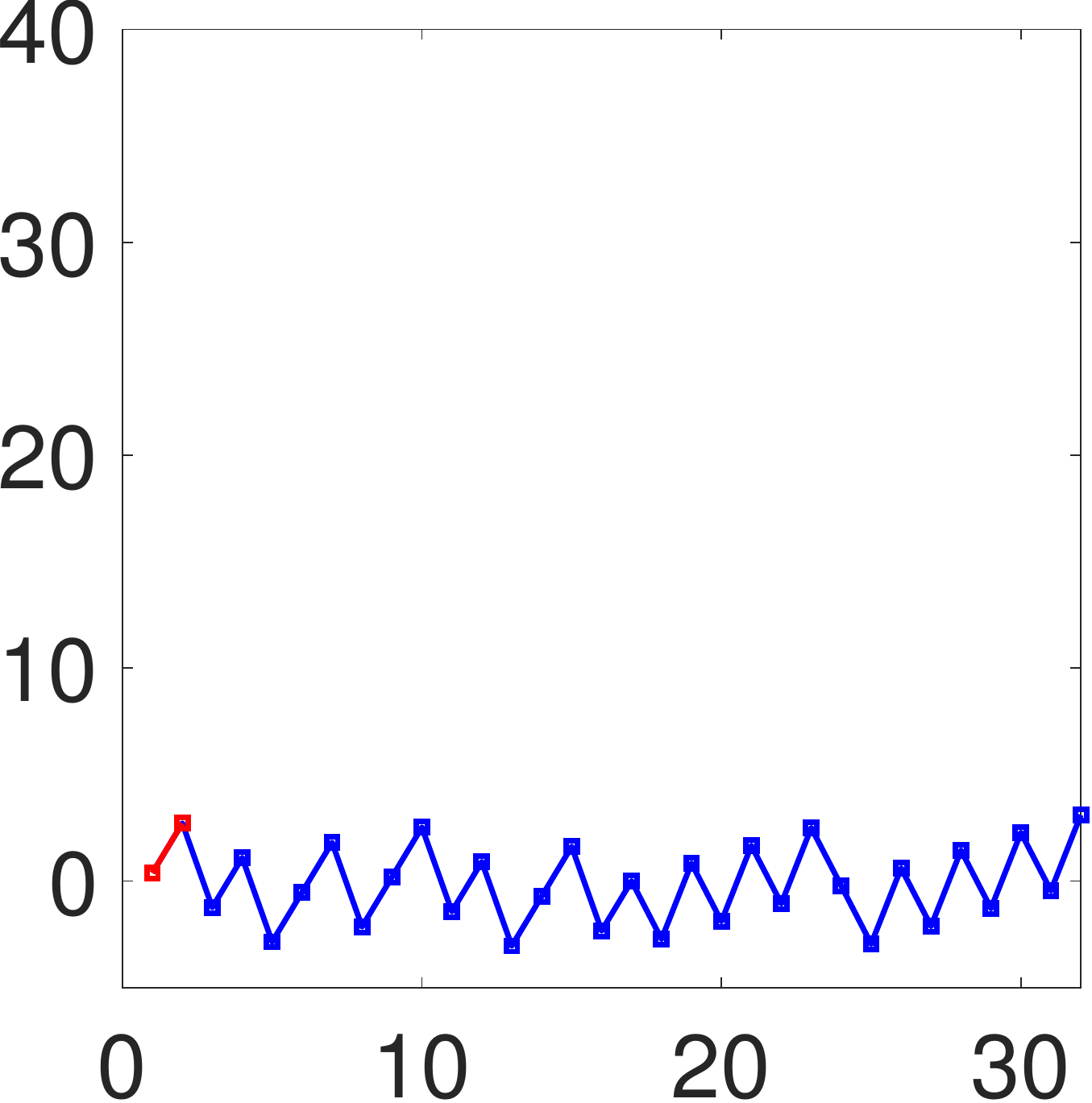}&\includegraphics[height=1in]{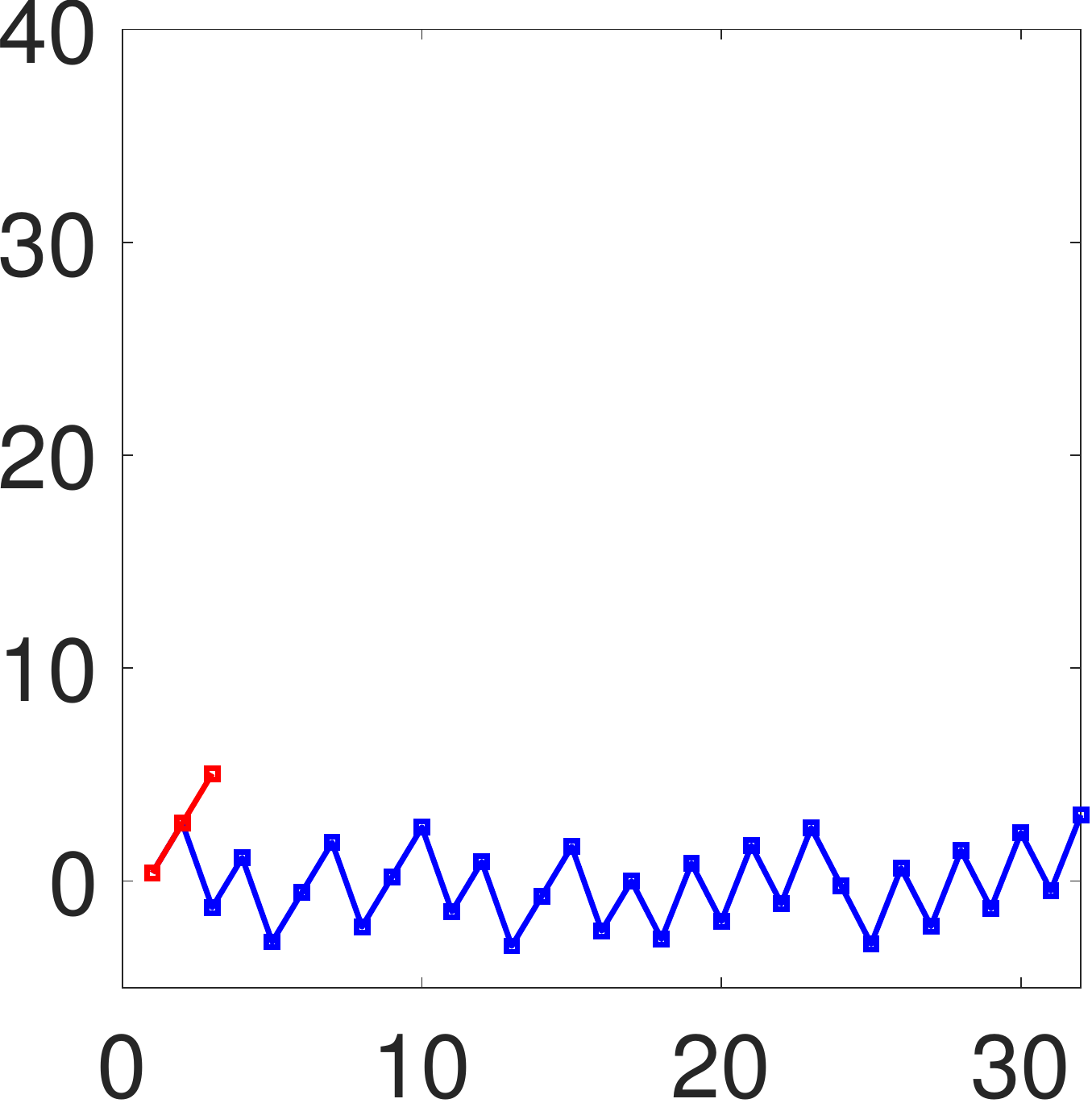}&\includegraphics[height=1in]{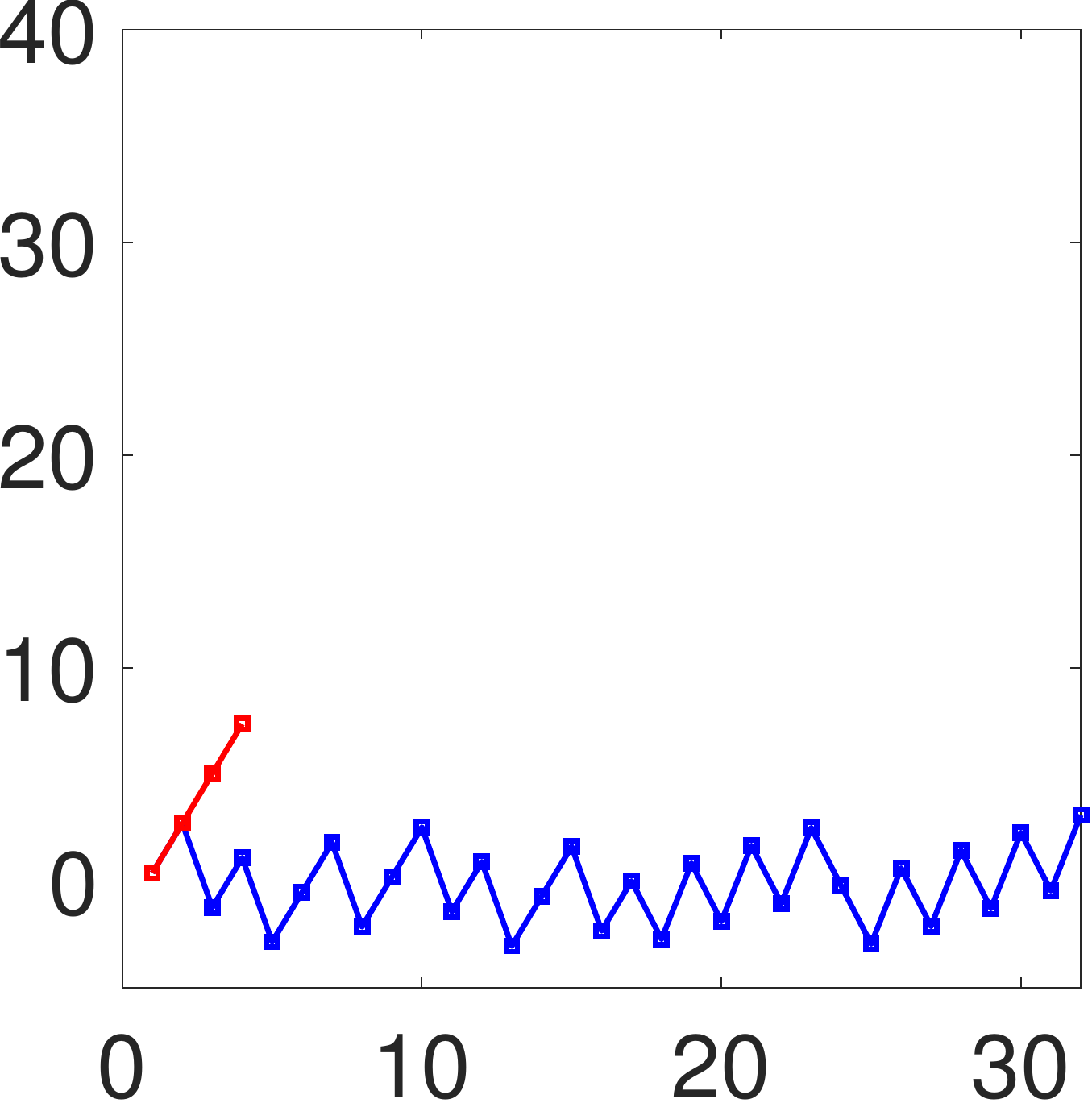}&\includegraphics[height=1in]{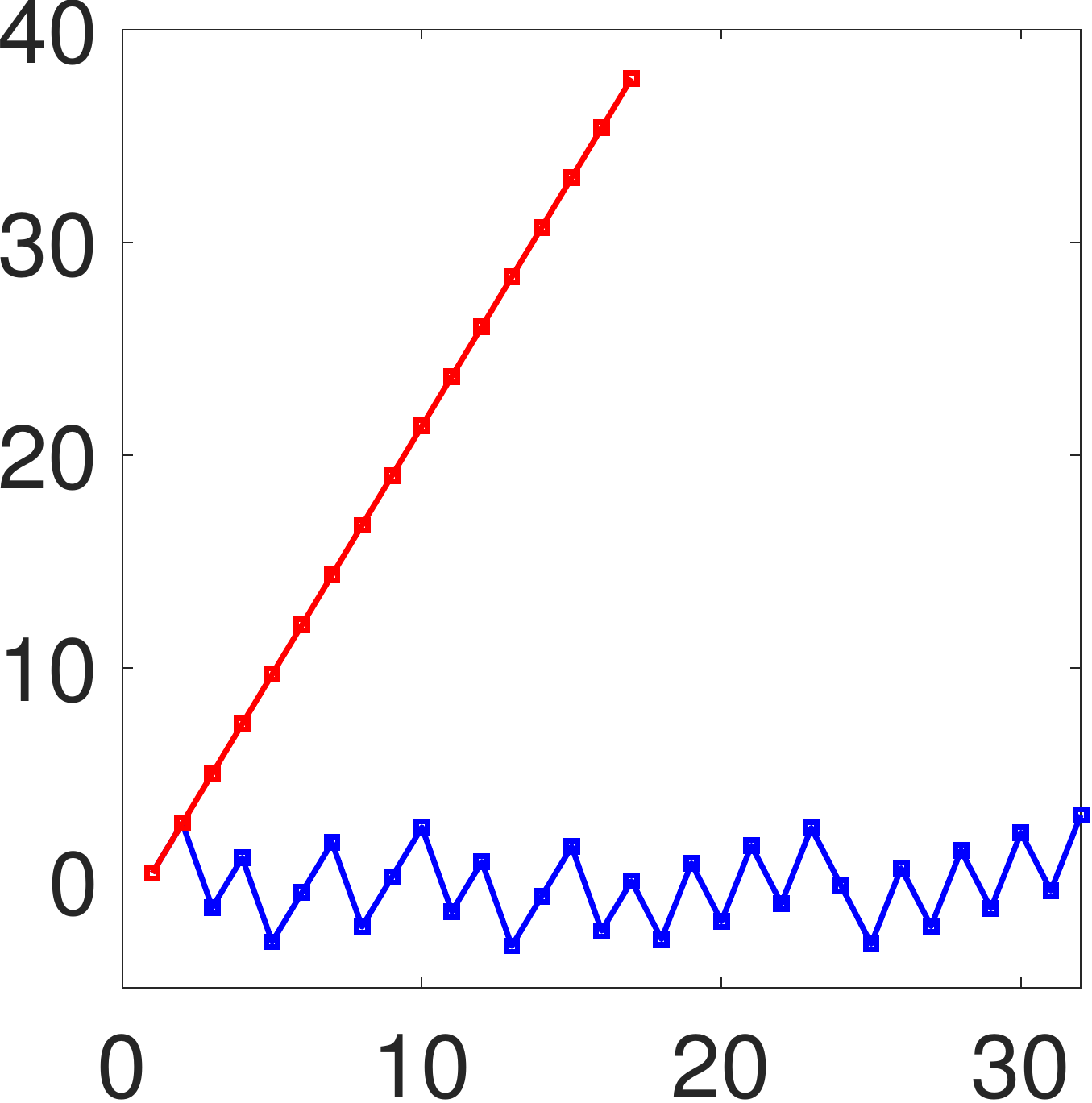}\\
      (a)&(b)&(c)&(d) &(e)\\
            \includegraphics[height=1in]{figure/p18.pdf}   & \includegraphics[height=1in]{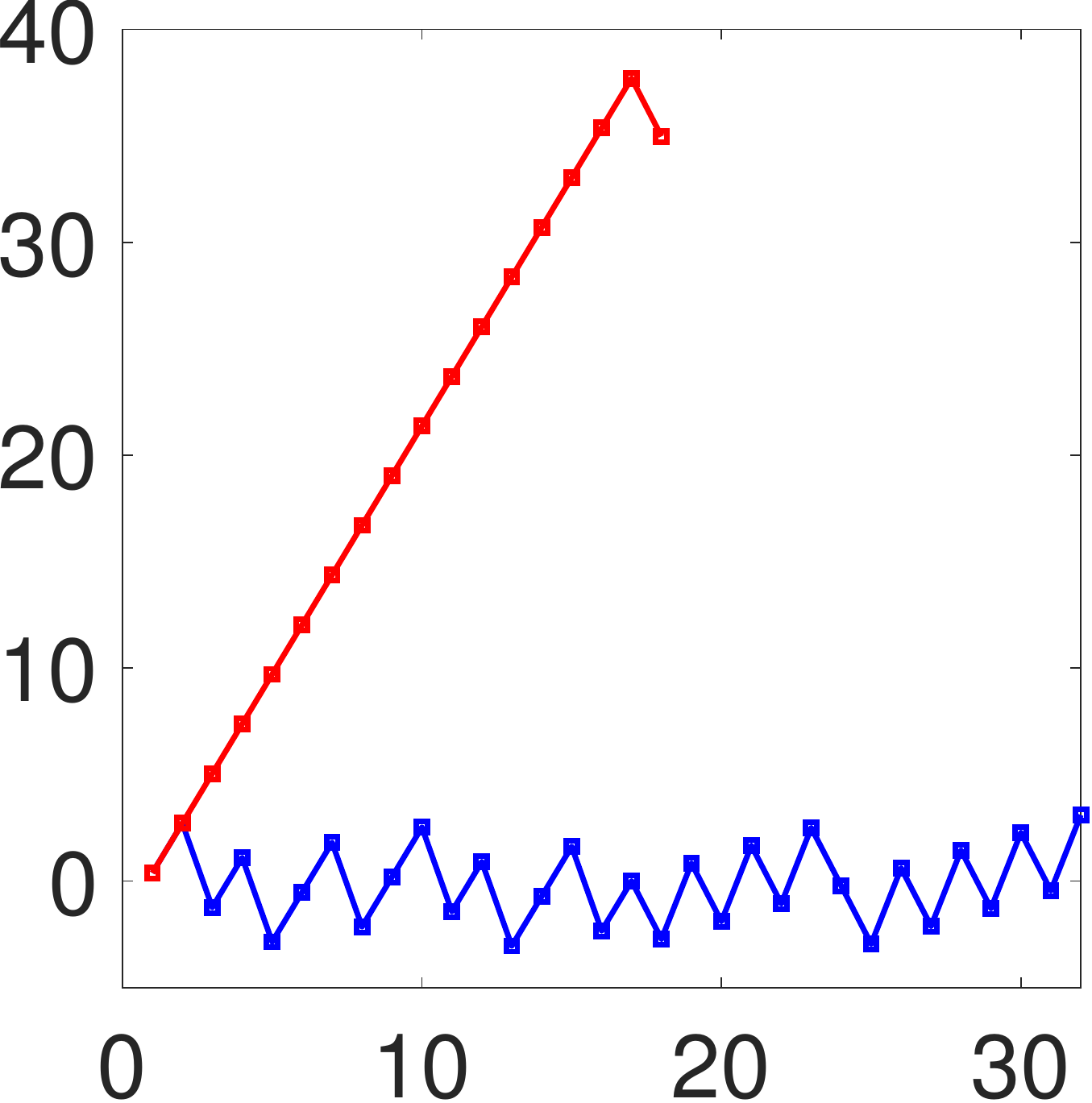}&\includegraphics[height=1in]{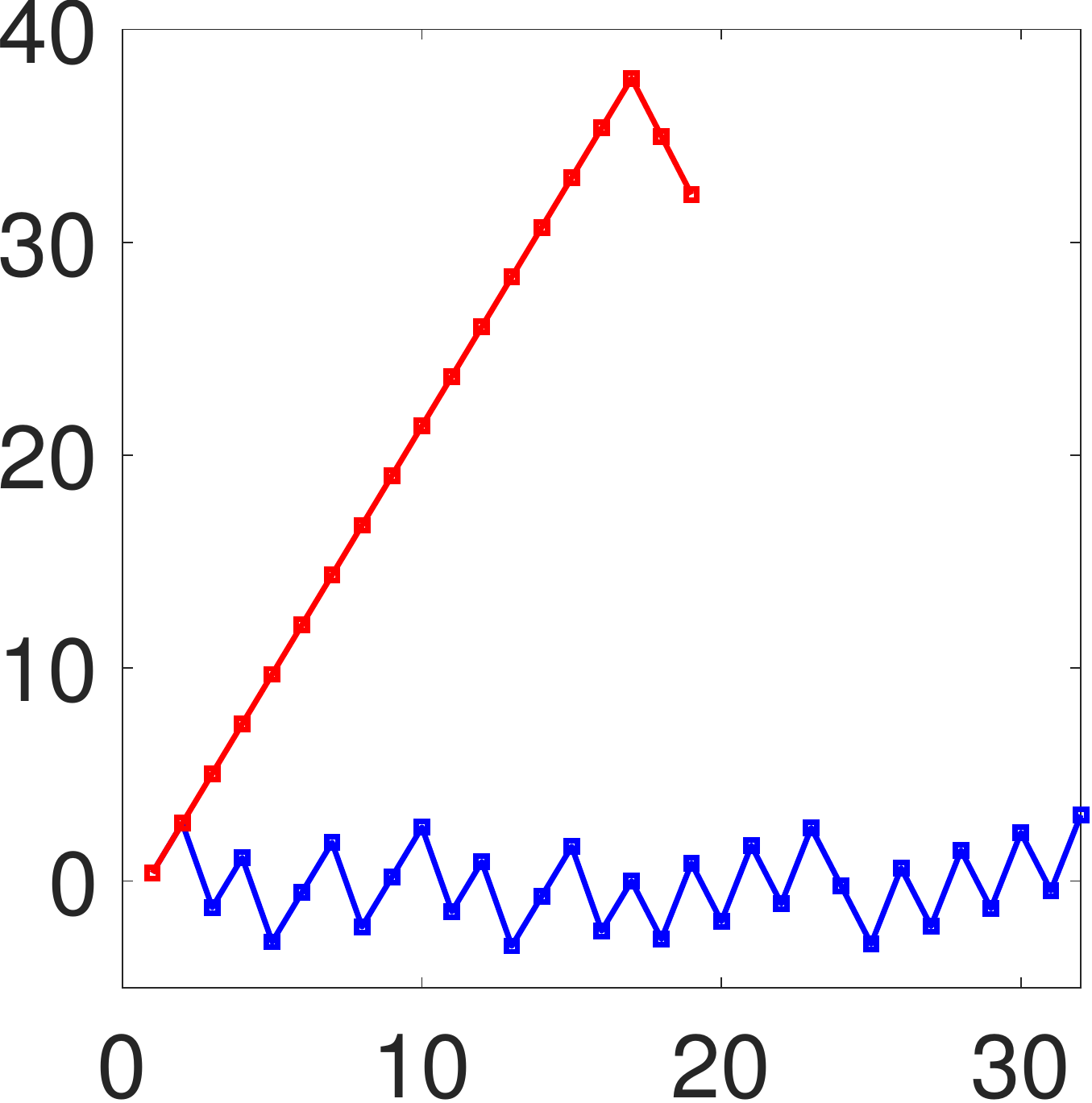}&\includegraphics[height=1in]{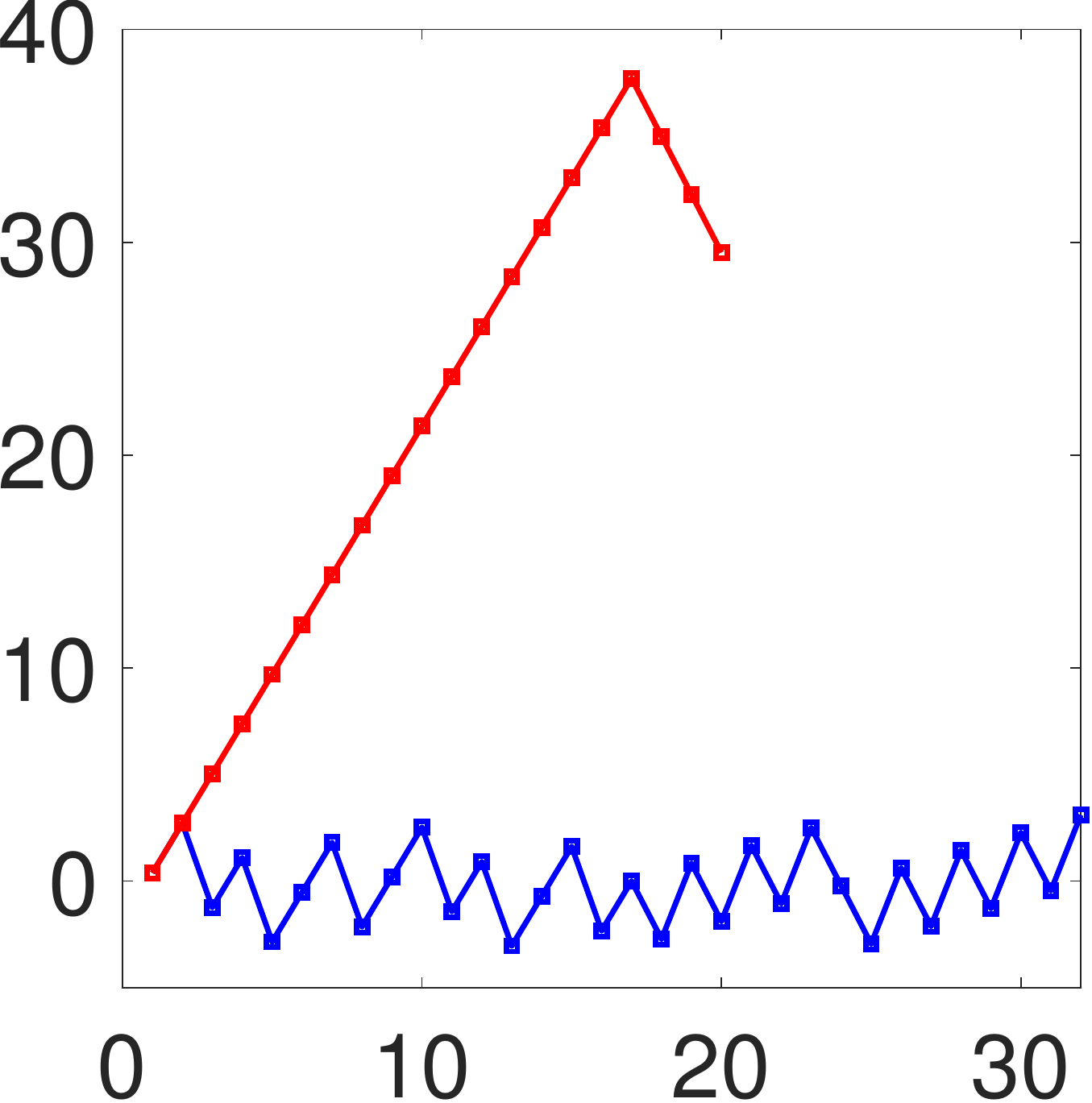}&\includegraphics[height=1in]{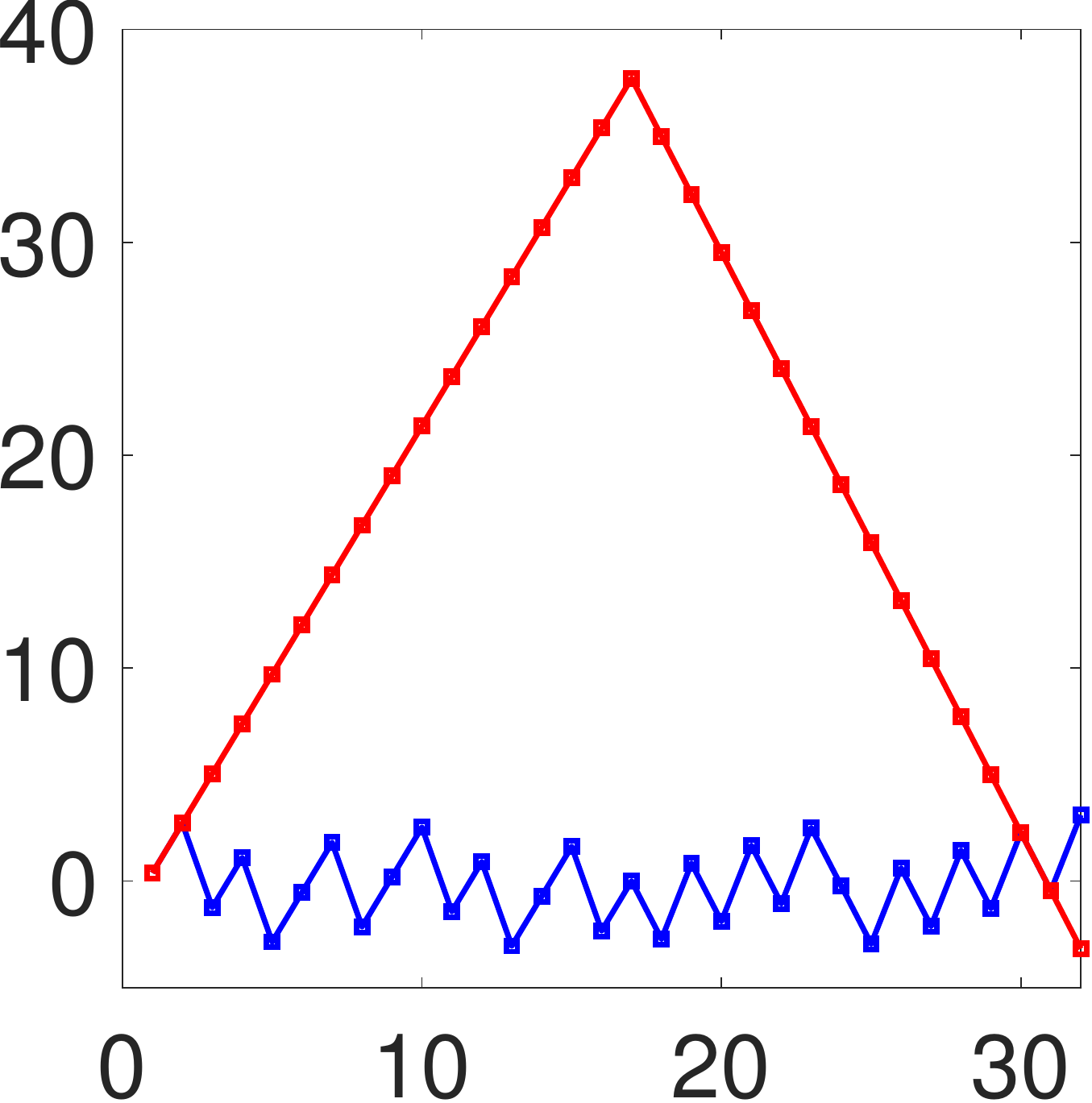}\\
            (f)&(g)&(h)&(i)&(j)
    \end{tabular}
  \end{center}
  \caption{Illustration of the recovery of one row of the phase function $\Phi(x,\xi) = x\cdot \xi + c(x)|\xi|$, where $c(x) = (2+\sin(2\pi x))/2$, by Algorithm \ref{alg:sst2}. This row is a function in $\xi$ denoted as $v$ of length $N$, and $v$ has two discontinuous point: one at the beginning and one in the middle. Suppose $u=\mod(v,1)$, we only know $u$ (in blue) and would like to recover $v$ (in red) from $u$. Top panel: (a) $u$.  (b) Line $6$-$9$ in Algorithm \ref{alg:sst2} assign the first two entries of $v$ right after the first discontinuous point such that they have the minimum distance while maintaining $\mod(v,1)=\mod(u,1)$. (c) Line $10$-$13$ in Algorithm \ref{alg:sst2} assign the third entry of $v$ such that $v(2)-v(1)$ and $v(3)-v(2)$ have the minimum distance while maintaining $\mod(v,1)=\mod(u,1)$.  (d) Line $16$ in Algorithm \ref{alg:sst2} assigns the fourth entry of $v$ such that $v(3)+v(1)-2v(2)$ and $v(4)+v(2)-2v(3)$ have the minimum distance while maintaining $\mod(v,1)=\mod(u,1)$. (e) Similarly, for all other $i$'s before the second discontinuous point, assign the $i$-th entry of $v$ by minimizing the distance between $v(i-1)+v(i-3)-2v(i-2)$ and $v(i)+v(i-2)-2v(i-1)$ while maintaining $\mod(v,1)=\mod(u,1)$.  Bottom panel: the second discontinuous point is detected by Line $17$ in Algorithm \ref{alg:sst2}; apply the same procedure as for (a)-(e) to recover the second part of $v$ after the second discontinuous point.}
\label{fig:vec}
\end{figure}

 A heuristic solution of the $TV^3$-norm minimization is to trace the columns and rows of $\frac{1}{2\pi} \Im\left(\log \left(K\right)\right)$ to identify smooth columns and rows of $\Psi$ agreeing with \eqref{eqn:tr} and satisfying the following conditions:
  \begin{enumerate}
  \item the variation of these columns and rows of $\Psi$ is small;
  \item recovered columns and rows after tracing share the same value at the intersection.
  \end{enumerate}
 {  Let us start with an example of vector recovery with $TV^3$-norm minimization to motivate the algorithm for matrix recovery:}
   \begin{eqnarray}
 \label{eqn:mv}
& \smash{\displaystyle\min_{v\in\mathbb{R}^{N}}} & \|v\|_{TV^3}\\
& \text{subject to} &  \mod(v,1) = \frac{1}{2\pi} \Im\left(\log \left(k\right)\right),\nonumber
 \end{eqnarray}
 {where $k\in\mathbb{R}^N$ is a given vector. The discussion below will be summarized in Algorithm \ref{alg:sst2}. Figure \ref{fig:vec} visualizes the vector recovery procedure for a simple case when $k$ is a vector from the discretization of $e^{2\pi i \phi(\xi)}$ with a piecewise smooth function $\phi(\xi)$ with only one discontinuous location $\xi=0$.}
 
{First, we assume $k$ is a vector from the discretization of $e^{2\pi i \phi(\xi)}$ with a smooth function $\phi(\xi)$. Let $u= \frac{1}{2\pi} \Im\left(\log \left(k\right)\right)$. We 
only know $u$ and would like to recover $v$ from $u$. If we have known $v(i:i+2)$, to minimize the $TV^3$-norm of $v$, we can assign the value of $v(i+3)$ such that $v(i+2)+v(i)-2v(i+1)$ and $v(i+3)+v(i+1)-2v(i+2)$ have the minimum distance while maintaining $\mod(v(i+3),1)=u(i+3)$ (corresponding to Line $16$ in Algorithm \ref{alg:sst2}). Hence, we only need to determine the values of $v(1:3)$ as the initial condition of the $TV^3$-norm vector recovery (corresponding to Line $6$-$13$ in Algorithm \ref{alg:sst2}). Similarly, to maximize the smoothness of $v$, we can assign the value of $v(i+2)$ such that $v(i+1)-v(i)$ and $v(i+2)-v(i+1)$ have the minimum distance while maintaining $\mod(v(i+2),1)=u(i+2)$ (corresponding to Line $10$-$13$ in Algorithm \ref{alg:sst2}). Finally, we can assign any value to $v(1)$ and determine the value of $v(2)$ such that $|v(2)-v(1)|$ is minimized with $\mod(v(2),1)=u(2)$ (corresponding to Line $6$-$9$ in Algorithm \ref{alg:sst2}). }

{Second, we deal with the case when $k$ is a vector from the discretization of $e^{2\pi i \phi(\xi)}$ with a piecewise smooth function $\phi(\xi)$. Suppose \[\mathcal{S}=\{c_1,c_2,\dots, c_n\}\] is an index set storing the discontinuity locations of $\phi(\xi)$ with $c_1=1<c_2<\dots<c_n<N$. $c_1=1$ since we can always assume that $\phi(\xi)$ is discontinuous at the end points of its domain. We can apply the algorithm just above to recover each piece $v(c_i:c_{i+1})$ for $i=1$, $\dots$, $n$. When $i=1$, we are free to set up any value for $v(c_1)$, while when $i>1$, $v(c_i)$ has been assigned according to the recovery for the previous piece corresponding to $v(c_{i-1})$. This difference is considered in the ``if" statement in Line $6$ and $10$ in Algorithm \ref{alg:sst2}. Since there is no prior information about $\mathcal{S}$ except that we know $c_1=1\in\mathcal{S}$, Algorithm \ref{alg:sst2} automatically determine the discontinuous locations in Line $17$-$19$ according to a threshold $\tau$: when the second derivative of $v$ at a certain location is larger than $\tau$, we consider $v$ is discontinuous at this location. }

{Recall the goal of matrix recovery in \eqref{eqn:tr}, it is not necessary to tune the parameter $\tau$ such that the discontinuous locations are exactly identified. If Algorithm \ref{alg:sst2} miss some discontinuous locations, Algorithm \ref{alg:sst} will provide a smoother estimation of the phase matrix; if Algorithm \ref{alg:sst2} artificially detects $O(1)$ fake discontinuous locations, Algorithm \ref{alg:sst} will provide an estimation of the phase matrix with more pieces of smooth domains. As long as \eqref{eqn:tr} is satisfied, all these estimations are satisfactory. In our numerical tests, $\tau$ is set to be $\frac{\pi}{2}$ for all numerical examples. Other values of $\tau$ result in similar numerical results, as long as $\tau$ is not close to $0$ such that there are too many fake discontinuous points that bring down the efficiency of Algorithm \ref{alg:sst}.}

\vspace{0.25cm}
\begin{algorithm2e}[H]
 \label{alg:sst2}

 \caption{An $O(N)$ algorithm for recovering a vector $v$ from the observation $u=\mod(v,1)$. The discontinuous locations of $v$ is automatically detected. See Figure \ref{fig:vec} for an illustration with a simple example.} 
 
Input: a vector $u$ of length $N$, a discontinuity detection parameter $\tau$.

Output: a vector $v$ satisfying $\mod(v,1)=\mod(u,1)$, and a vector of indices $\mathcal{S}$ for discontinuity locations.

Initialize: $\mathcal{S}=[1]$; let $n$ be the number of elements in $\mathcal{S}$; and let $c=1$.

\While{$c\leq n$}{

If $c<n$, let $st=\mathcal{S}(c)$ and $ed=\mathcal{S}(c+1)-1$; otherwise, let $st=\mathcal{S}(c)$ and $ed=N$.

\If{$c=1$}{
Assign the values of $v(st:st+1)$ such that these two values have the minimum distance while maintaining $\mod(u(st:st+1),1)=\mod(v(st:st+1),1)$.
}
\Else{
Assign the values of $v(st)$ such that two values in $v(st-1:st)$ have the minimum distance while maintaining $\mod(u(st-1:st),1)=\mod(v(st-1:st),1)$.
}

\If{$c=1$}{
Assign the values of $v(st+2)$ such that $v(st+2)-v(st+1)$ and $v(st+1)-v(st)$ have the minimum distance while maintaining $\mod(u(st+2),1)=\mod(v(st+2),1)$.
}
\Else{
Assign the values of $v(st+1)$ such that $v(st+1)-v(st)$ and $v(st)-v(st-1)$ have the minimum distance while maintaining $\mod(u(st+1),1)=\mod(v(st+1),1)$.
}

If $c=1$, let $bg=st+3$; otherwise, let $bg=st+2$.

\For{all indices $a$ from $bg$ to $ed$}{
Assign the value of $v(a)$ such that $v(a-1)+v(a-3)-2v(a-2)$ and $v(a)+v(a-2)-2v(a-1)$ have the minimum distance while maintaining $\mod(v(a),1)=\mod(u(a),1)$. 

\If{$|v(a)+v(a-2)-2v(a-1)|>\tau$}{
Consider $a$ as a new location at which $v$ is discontinuous, add $a$ to $\mathcal{S}$, and let $n\leftarrow n+1$.

Break the for-loop.
}
}

{$c\leftarrow c+1$}.
}
\end{algorithm2e}
  \vspace{0.25cm}
  
 { When the vector recovery algorithm in Algorithm \ref{alg:sst2} is ready, we apply it to design a matrix recovery algorithm in Algorithm \ref{alg:sst}. Recall that recovered columns and rows by Algorithm \ref{alg:sst2} should share the same value at the intersection. To guarantee this, we carefully choose the recovery order of the rows and columns, and the initial values of vector recovery, to avoid assignment conflicts at the intersection. For simplicity, we only introduce Algorithm \ref{alg:sst} for a phase function defined on $\mathbb{R}\times \mathbb{R}$. We will leave the extension to high-dimensional case as a future work.}

\vspace{0.25cm}
\begin{algorithm2e}[H]
 \label{alg:sst}

 \caption{An $O(N)$ algorithm for the approximate solution of the $TV^3$-norm minimization when the phase function $\Phi(x,\xi)$ is defined on $\mathbb{R}\times \mathbb{R}$.} 
 
Input: a vector $\mathcal{C}$ and a vector $\mathcal{R}$ as the column and row index sets indicating $O(1)$ randomly selected columns and rows of $\Phi$, columns $U=\mod(\Phi(:,\mathcal{C}),1)$, rows $V=\mod(\Phi(\mathcal{R},:),1)$, a discontinuity detection parameter $\tau$.

Output: columns $\bar{U}$ and rows $\bar{V}$ satisfying $\mod(\bar{U},1)=\mod(U,1)$, and $\mod(\bar{V},1)=\mod(V,1)$.

Apply Algorithm \ref{alg:sst2} to $U(:,\mathcal{C}(1))$ to detect a discontinuous point set $\mathcal{S}_r$; add $\mathcal{S}_r$ to $\mathcal{R}$ and update row samples $V$ accordingly.

Apply Algorithm \ref{alg:sst2} to $V(\mathcal{R}(1),:)$ to detect a discontinuous point set $\mathcal{S}_c$; add $\mathcal{S}_c$ to $\mathcal{C}$ and update row samples $U$ accordingly.

Let $n_r$ be the number of elements in $\mathcal{S}_r$ and $n_c$ be the number of elements in $\mathcal{S}_c$. The discontinuous point sets naturally partition the phase matrix into $n_r\times n_c$ blocks (see Figure \ref{fig:trace} for an example).

\For{Each block partitioned by discontinuous point sets}{
Set $\tau=2\pi$, since it is not necessary to detect discontinuity here.

Apply Algorithm \ref{alg:sst2} to recover the first row and the first column of each block.

Apply Algorithm \ref{alg:sst2} to recover the second and the third columns of each block. Make sure that the recovery shares the same entries when they intersect with the first row, and there is no discontinuity along rows inside the first three columns. 

Apply Algorithm \ref{alg:sst2} to recover $O(1)$ rows of each block such that the first three entries of these rows have the same entries as in the first three columns. 

Apply Algorithm \ref{alg:sst2} to recover $O(1)$ columns of each block such that these columns have the same entries as in the recovered rows when a column and a row intersects.
}

\end{algorithm2e}
\vspace{0.25cm}

 \begin{figure}[ht!]
  \begin{center}
    \begin{tabular}{cccc}
      \includegraphics[height=1.2in]{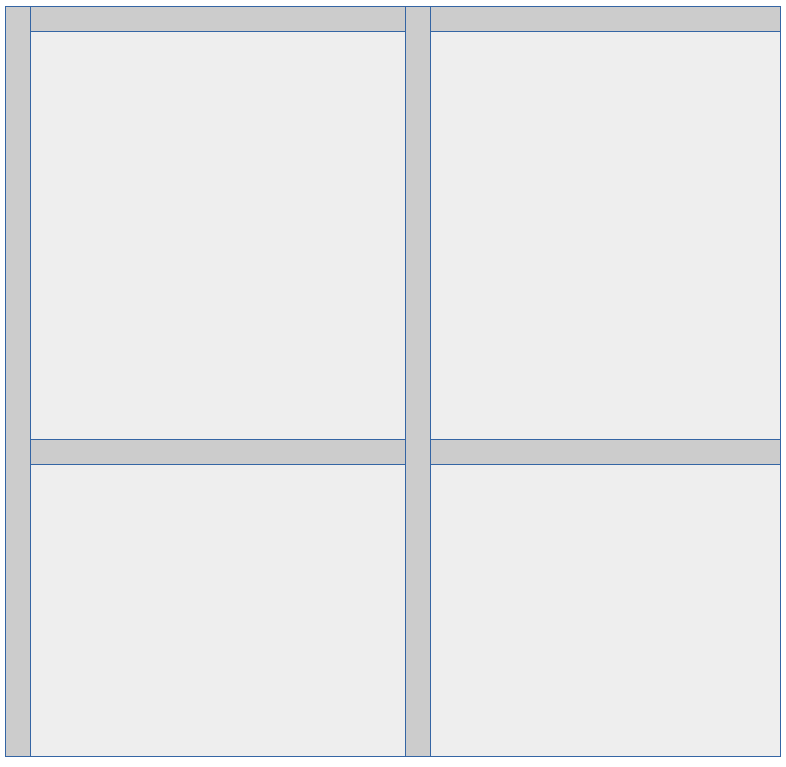} &\includegraphics[height=1.2in]{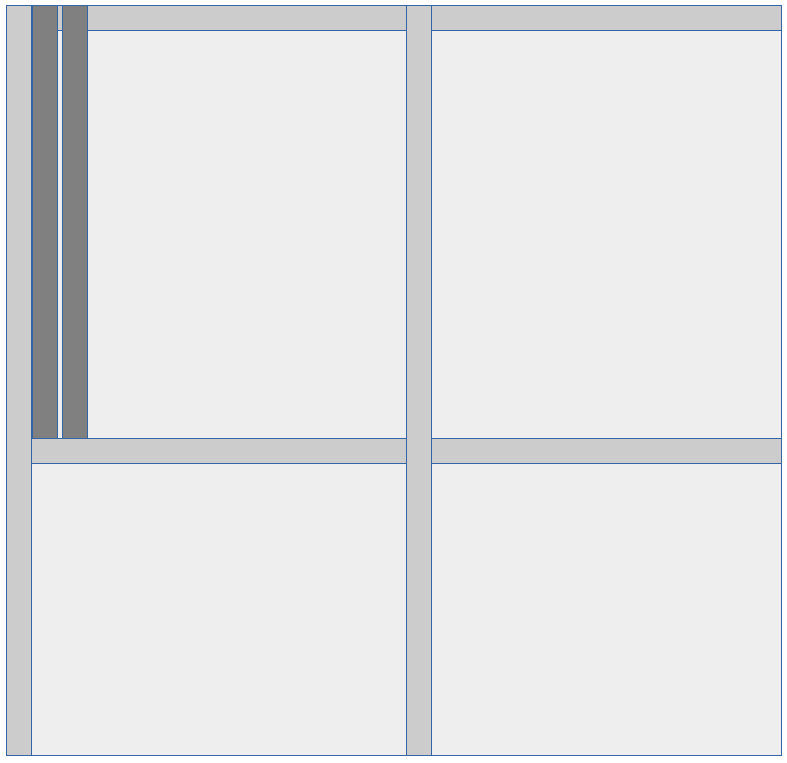}&   \includegraphics[height=1.2in]{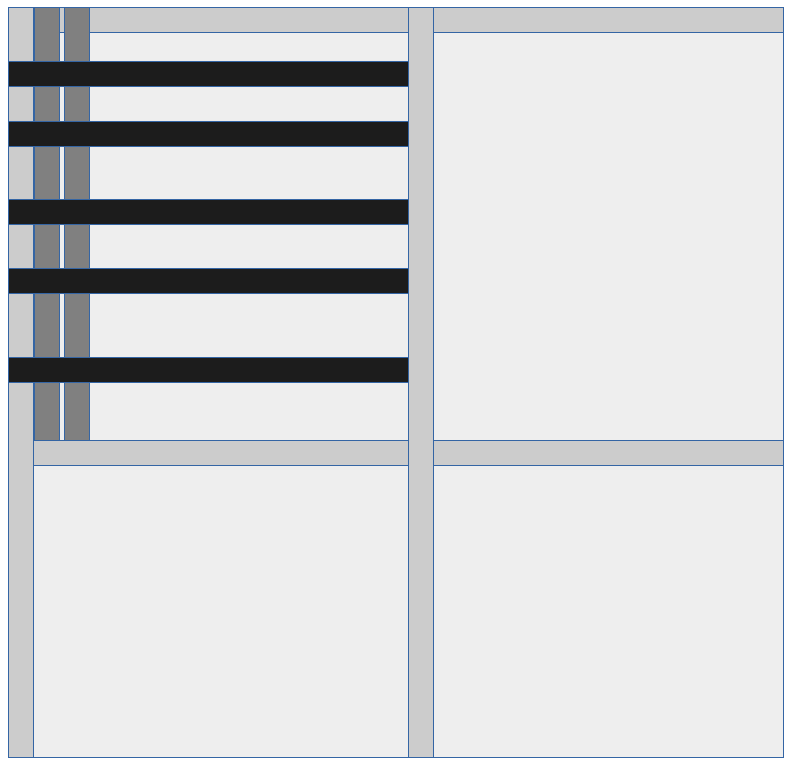} &\includegraphics[height=1.2in]{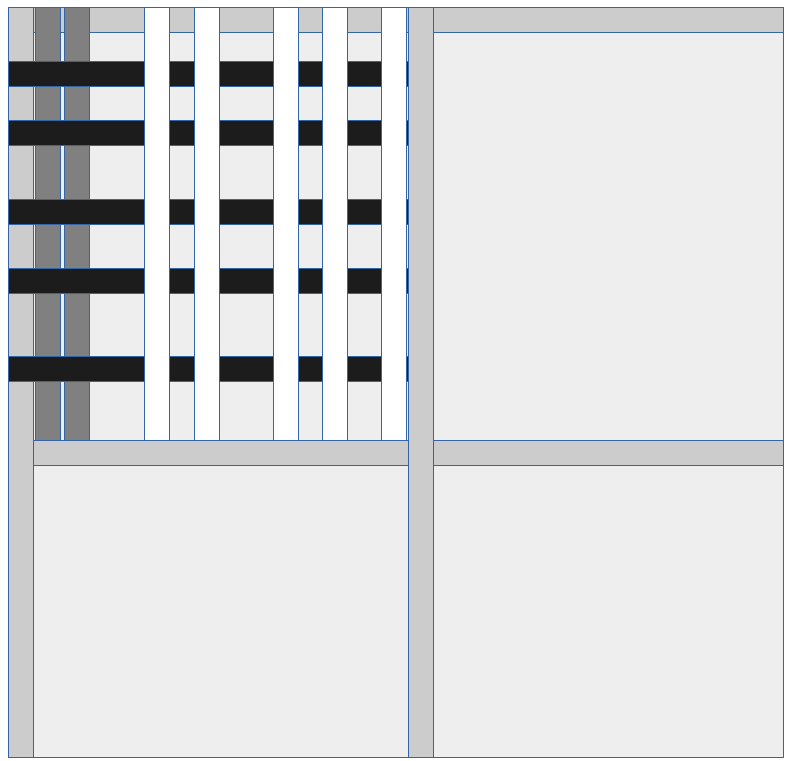}\\
      (a) & (b) & (c) & (d)
    \end{tabular}
  \end{center}
  \caption{Illustration of the low-rank matrix recovery in Algorithm \ref{alg:sst}. (a) The matrix is partitioned into submatrices such that there is no discontinuity along columns and rows in each submatrix. Line $8$ in Algorithm \ref{alg:sst} recovers the first column and row of each submatrix. (b) Next, Line $9$ in Algorithm \ref{alg:sst} recovers the second and the third columns for each submatrix. (c) Next, Line $10$ in Algorithm \ref{alg:sst} recovers $O(1)$ rows of each submatrix such that the first three entries of these rows have the same entries as in the first three columns. (d) Finally, Line $11$ in Algorithm \ref{alg:sst} recovers $O(1)$ columns of each submatrix such that these columns have the same entries as in the recovered rows when a column and a row intersects.}
\label{fig:trace}
\end{figure}

In the case of higher dimensions, the discretization of the oscillatory integral transform and the arrangement of grid points will lead to artificial discontinuity along the column and row indices. For example, a column or a row as a one-dimensional {function} in index is discontinuous at a certain point, while we look back to the original high dimensional domain, the original kernel function is continuous at the corresponding point. Hence, once the {discretization} and arrangement of grid points have been fixed, we can remove the artificial discontinuity and apply the same ideas as in Algorithm \ref{alg:sst} to recover high dimensional phase functions.  

 With Algorithm \ref{alg:sst} ready, we are able to introduce the nearly linear scaling algorithm for constructing a low-rank factorization $UV^*$, where $U\in\mathbb{C}^{N\times r}$ and $V\in\mathbb{C}^{N\times r}$,  such that $e^{2\pi\i  UV^*}=e^{2\pi\i  \Phi}$ when we only know the kernel matrix $K=e^{2\pi\i \Phi}$ through Scenarios $1$ and $2$ in Table \ref{tab:sc}. This method is summarized in Algorithm \ref{alg:idrLR}.

\subsection{Summary for the low-rank matrix factorization in the unified framework}

Before moving to the algorithms for other main steps of the unified framework as shown in Figure \ref{fig:flow}, let us summarize how those algorithms in Section \ref{sec:LRF} and Section \ref{sec:nLR} can be applied to construct the low-rank matrix factorization of the ampltiude and phase functions with nearly linear computational complexity.

For a general kernel $K(x,\xi)=\alpha(x,\xi)e^{2\pi\i \Phi(x,\xi)}$, suppose we discretize $\alpha(x,\xi)$ and $\Phi(x,\xi)$ with $N$ grid points in each variable to obtain the amplitude matrix $\mathcal{A}$ and the phase matrix $\Phi$. When the explicit formulas of $\alpha(x,\xi)$ and $\Phi(x,\xi)$ are known, it takes $O(N)$ operations to evaluate one column or one row of $\mathcal{A}$ and $\Phi$. Hence, Algorithm \ref{alg:rlr} in Section \ref{sec:LRF} is able to construct the low-rank matrix factorization of $\mathcal{A}$ and $\Phi$ in $O(N)$ operations. 

When the explicit formulas are unknown but they are solutions of certain PDE's as in Scenario $3$ in Table \ref{tab:sc}. In this paper, we simply assume that $O(1)$ columns and rows of the amplitude and phase functions are available and Algorithm \ref{alg:rlr} in Section \ref{sec:LRF} can be applied to construct the low-rank factorization in $O(N)$ operations. In practical applications like solving wave equations \cite{Yingwave}, this assumption for the phase function is reasonable since it can be obtain via interpolating the solution of the PDE's on a coarse grid of size independent of $N$. However, obtaining the amplitude function might take expensive computation for solving PDE's on a grid depending on $N$. Optimizing this complexity will be left as interesting future work.

In the case of indirect access in Scenario $1$ and $2$ in Table \ref{tab:sc}, it takes $O(N)$ or $O(N\log N)$ operations to evaluate one column or one row of the kernel matrix $K$. By taking the absolute value of $K$, we obtain one column or one row of $\mathcal{A}$. Hence, the low-rank factorization of $\mathcal{A}$ can be constructed via Algorithm \ref{alg:rlr} in Section \ref{sec:LRF} in $O(N\log N)$ operations. Dividing the amplitude from the kernel, we have the access of the phase in the form of $e^{2\pi\i  \Phi(x,\xi)}$. Hence, the low-rank factorization of $\Phi$ can be constructed by Algorithm \ref{alg:idrLR} in Section \ref{sec:nLR} in $O(N\log N)$ operations.

\vspace{0.5cm}
\begin{algorithm2e}[H]
 \label{alg:idrLR}
 Input: Scenario $1$: an algorithm for evaluating an arbitrary entry of the kernel matrix $K$ in $O(1)$ operations; Scenario $2$: an $O(N\log N)$ algorithm for applying $K$ and its transpose to a vector. A rank parameter $r$, an over-sampling parameter $q$, and the matrix size $N$.

Output: $U\in\mathbb{C}^{N\times r}$ and $V\in\mathbb{C}^{N\times r}$ such that $e^{2\pi\i UV^*}=e^{2\pi\i  \Phi}$.

\caption{Low-rank matrix factorization for indirect access. The computational complexity in Scenario $1$ is $O(N)$ and that in Scenario $2$ is $O(N\log N)$.}
 \uIf{Scenario $1$}{Evaluate $rq$ randomly selected columns and rows of $K$.}
 \uElseIf{Scenario $2$}{Apply the kernel matrix $K$ and its transpose to $rq$ randomly selected natural basis vectors in $\mathbb{R}^N$ to obtain the columns and rows of $K$.}
 
 Apply Algorithm \ref{alg:sst} with the columns and rows of $K$ to obtain $rq$ columns and rows of a matrix $\Psi$ such that $e^{2\pi\i  \Psi}=e^{2\pi\i  \Phi}$.
 
 Apply Algorithm \ref{alg:rlr} with the columns and rows of $\Psi$ to obtain the low-rank factorization of $\Psi\approx UV^*$ such that $e^{2\pi\i UV^*}=e^{2\pi\i  \Phi}$, $U\in\mathbb{C}^{N\times r}$, and $V\in\mathbb{C}^{N\times r}$.
 
\end{algorithm2e}

\vspace{0.25cm}

\section{NUFFT and dimension lifting}
\label{sec:frame}
This section introduces a new NUFFT approach by dimension lifting to evaluate the oscillatory integral transform 
\begin{equation}
\label{eqn:tg2}
g(x) = \int \alpha(x,\xi)e^{2\pi\i \Phi(x,\xi)} f(\xi)d\xi.
\end{equation}
If we could find $\{p_j(x)\}_{1\leq j\leq r}$ and $\{q_j(\xi)\}_{1\leq j\leq r}$ such that $e^{2\pi\i (\Phi(x,\xi)-\sum_{j=1}^r p_j(x)q_j(\xi))}$ is numerically low-rank, then \eqref{eqn:tg2} is reduced to $O(1)$ $r$-dimensional NUFFT's:
\begin{equation}
\label{eqn:tg3}
g(x)\approx \sum_{k=1}^{r_\epsilon} a_k(x)\int e^{2\pi\i  \sum_{j=1}^r p_j(x)q_j(\xi)}\left(b_k(\xi)f(\xi)\right)d\xi,
\end{equation}
where $a_k(x)$ and $b_k(\xi)$ are the low-rank approximation of 
\[
\alpha(x,\xi)e^{2\pi\i (\Phi(x,\xi)-\sum_{j=1}^r p_j(x)q_j(\xi))}\approx \sum_{k=1}^{r_\epsilon} a_k(x)b_k(\xi).
\]
Note that the prefactor of an $r$-dimensional NUFFT increases as $r$ increases. {Hence, the key condition for deciding whether NUFFT is suitable for evaluating \eqref{eqn:tg2} is the existence of $\{p_j(x)\}_{1\leq j\leq r}$ and $\{q_j(\xi)\}_{1\leq j\leq r}$ to ensure a small $r$ and $r_\epsilon$.}

The choice of $\{p_j(x)\}_{1\leq j\leq r}$ and $\{q_j(\xi)\}_{1\leq j\leq r}$ is related to but different from classical low-rank approximation problems that can be solved by the SVD. In fact, we have a new low-rank approximation problem for fixed rank parameters $r$ and $r_\epsilon$ as follows:
 \begin{eqnarray}
 \label{eqn:minexp}
& \smash{\displaystyle\min_{P,Q\in\mathbb{R}^{N\times r},U,V\in\mathbb{R}^{N\times r_\epsilon}}} & \|\mathcal{A}e^{2\pi\i (\Phi-PQ^*)}-UV^*\|_2,
 \end{eqnarray} 
 where $\mathcal{A}$ represents the amplitude matrix for $\alpha(x,\xi)$, and $\Phi$ is the phase matrix for $\Phi(x,\xi)$. An immediate idea is to set reasonable $r$ and $r_\epsilon$, and solve the minimization problem in \eqref{eqn:minexp}. If the minimum value of the objective function is sufficiently small, then we can use the NUFFT to evaluate \eqref{eqn:tg2} via \eqref{eqn:tg3}. However, solving the optimization problem in \eqref{eqn:minexp} could be much more expensive than $O(N)$. This motivates Algorithm \ref{alg:dc} below for deciding whether we could use NUFFT in $O(N)$ operations.

\vspace{0.5cm}
\begin{algorithm2e}[H]
 \label{alg:dc}
 Input: low-rank factorization of the phase matrix $\Phi\approx U_1V_1^*$, where $U_1\in\mathbb{C}^{N\times r_1}$ and $V_1\in\mathbb{C}^{N\times r_1}$; low-rank factorization of the amplitude matrix $\mathcal{A}\approx U_2V_2^*$, where $U_2\in\mathbb{C}^{N\times r_2}$ and $V_2\in\mathbb{C}^{N\times r_2}$; rank parameters $r<r_1$ and $r_\epsilon$, an over-sampling parameter $q>1$, an accuracy parameter $\epsilon\approx 0$, and the matrix size $N$.

Output: $y\in\{0,1\}$; if $y=1$, return $P,Q\in\mathbb{R}^{N\times r},U,V\in\mathbb{R}^{N\times r_\epsilon}$ satisfying the low-rank factorization \[(U_2V_2^*)\odot e^{2\pi\i (U_1V_1^*-P Q^*)}\approx UV^*,\]
where $\odot$ means the entry-wise dot product of two matrices.

\caption{An $O(N)$ algorithm for deciding whether NUFFT is applicable; if NUFFT is applicable, returns the low-rank factorization for the evaluation in \eqref{eqn:tg3}.}

Compute the approximate $r$-leading SVD of the rank-$r_1$ matrix $U_1V_1^*$ using the randomized truncated SVD algorithm in \cite{Gu,randomSVD}\footnote{In the computation of the leading singular pair, since we have the rank-$r_1$ factorization, the computational cost is $O(N)$, the convergence to the ground true singular pair is very fast if a test matrix with a number of columns larger than $r_1$ is applied \cite{Gu}, and the probability to obtain high accuracy is very close to $1$.} and denote it as $P\Sigma Q^*$. Update $P\Sigma\rightarrow P$.

Randomly sample $rq$ columns of $(U_2V_2^*)\odot e^{2\pi\i (U_1V_1^*-P Q^*)}$ and stack them into a matrix $M$. Perform a pivoted {QR} decomposition of $M$ and let $R$ be the resulting $rq\times rq$ upper triangular matrix. 
  
Let $n$ be the number of diagonal entries of $R$ that are larger than $R(1,1)\epsilon$. If $n<r$, let $y=1$; otherwise, let $y=0$.

\If{$y=1$}{Apply Algorithm \ref{alg:rlr} to compute the low-rank factorization $UV^*$ of $(U_2V_2^*)\odot e^{2\pi\i (U_1V_1^*-P Q^*)}$ with the rank parameter $r_\epsilon$ and the over-sampling parameter $q$.}
\end{algorithm2e}

\vspace{0.5cm}

Although Algorithm \ref{alg:dc} is not optimal in the sense that it cannot provide the best $P$ and $Q$ such that the low-rank approximation of $\mathcal{A}e^{2\pi\i (\Phi-PQ^*)}$ has the smallest rank, Algorithm \ref{alg:dc} is sufficiently efficient and works well in practice. The stability and probability analysis of the main components of this algorithm can be found in \cite{Gu,randomSVD,RNL}. If the output of Algorithm \ref{alg:dc} is $y=1$, then the low-rank factorization of $\mathcal{A}e^{2\pi\i (\Phi-PQ^*)}$ is incorporated into \eqref{eqn:tg3} to evaluate \eqref{eqn:tg2} with $r_\epsilon$ $r$-dimensional NUFFT's. {Note that $r$ is a parameter less than or equal to $3$ according to the current development of NUFFT, and $r_\epsilon$ usually can be as large as $O(100)$ since $N$ is usually very large.} If the output of Algorithm \ref{alg:dc} is $y=0$, then the NUFFT approach is not applicable and we use the IBF-MAT introduced below to evaluate \eqref{eqn:tg2}. As we shall see later in the numerical examples, in some applications, the NUFFT approach is not applicable for the whole matrix $K$, but it can be applied to submatrices of $K$. Combining the results of all the submatrices of $K$ can also lead to efficient matvec for $K$. This strategy is problem-dependent and hence we omit the detailed discussion here.
 
\section{IBF-MAT}
\label{sec:BF}
This section introduces the IBF-MAT for evaluating the oscillatory integral transform if NUFFT is not applicable. Recall that after computing the low-rank factorization of the amplitude function, our target is to evalutate \eqref{eqn:tg}. If NUFFT is not applicable, we compute the IBF-MAT of $e^{2\pi\i  \Phi(x,\xi)}$, where the phase function is given in a form of a low-rank matrix factorization. Then the evaluation of \eqref{eqn:tg} is reduced to the application of IBF-MAT to $O(1)$ vectors. Hence, we only focus on the IBF-MAT of $e^{2\pi\i  UV^*}$, where $U$ and $V\in\mathbb{R}^{N\times r}$. To simplify the discussion, we also assume that $x$ and $\xi$ are one-dimensional variables. It is easy to extend the IBF-MAT to multi-dimensional cases following the ideas in \cite{FIO09,IBF,FIO14,MBF}.

\begin{figure}[ht!]
  \begin{center}
    \begin{tabular}{ccccc}
      \includegraphics[height=1.1in]{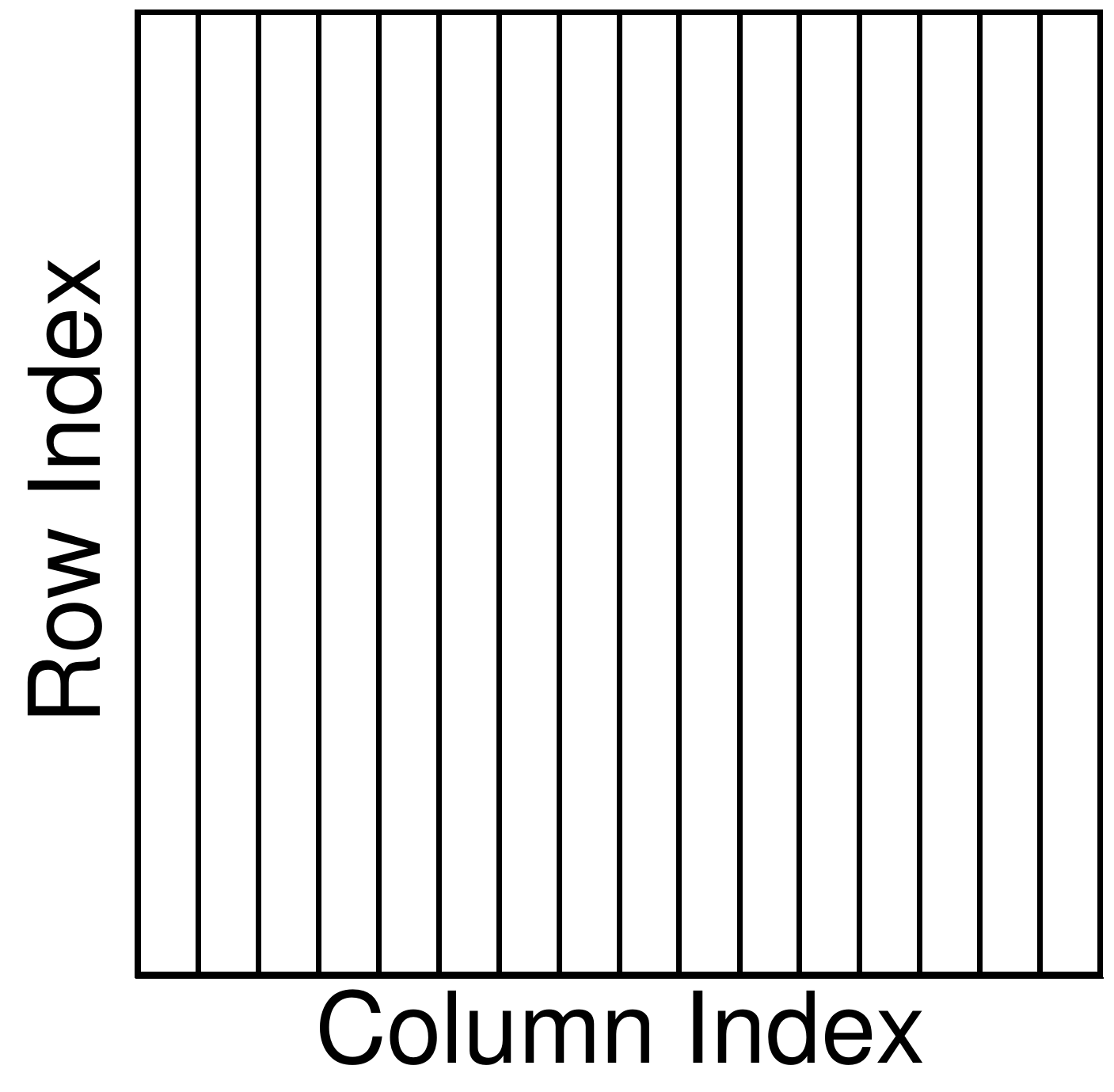}&
      \includegraphics[height=1.1in]{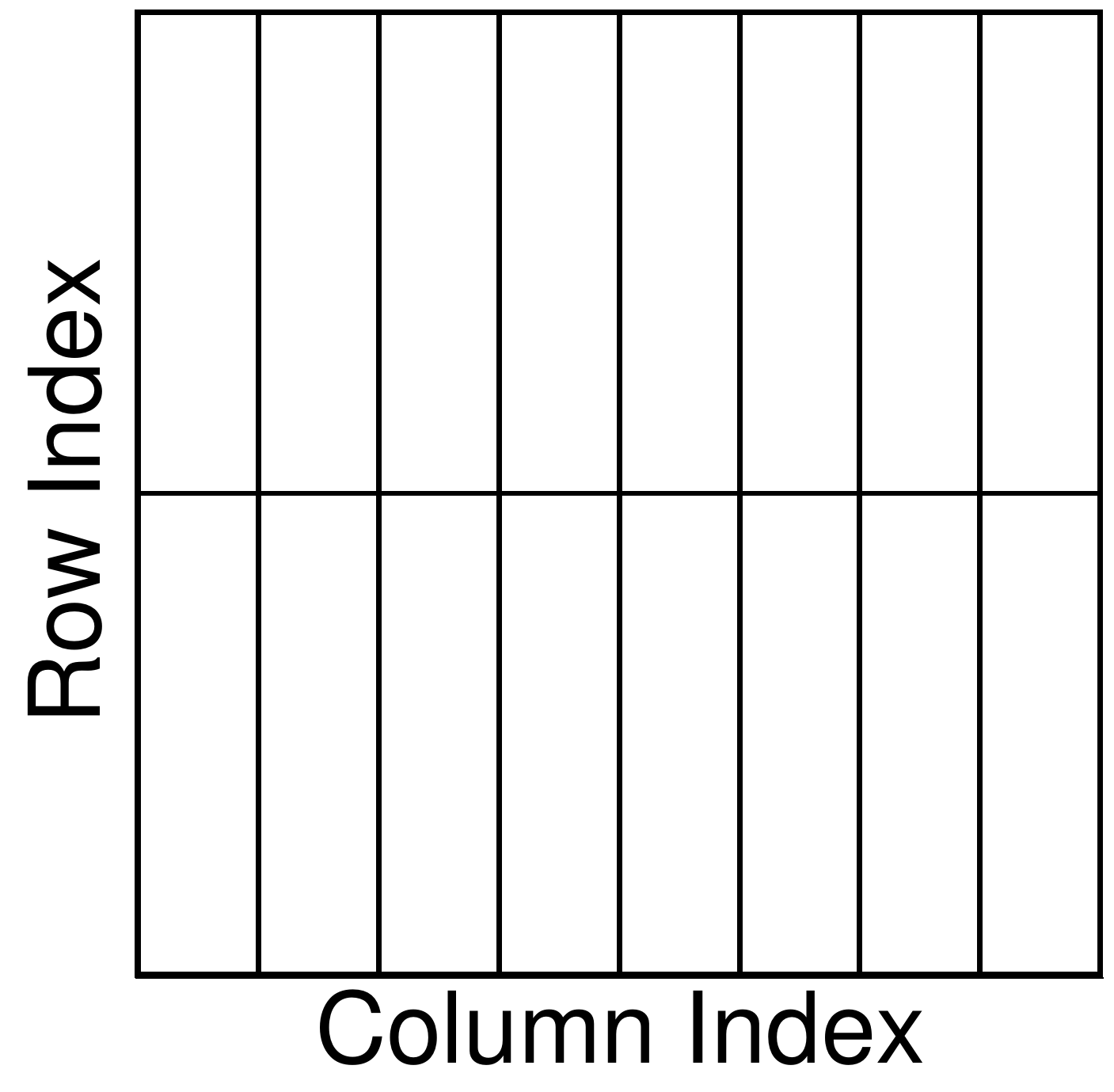}&
      \includegraphics[height=1.1in]{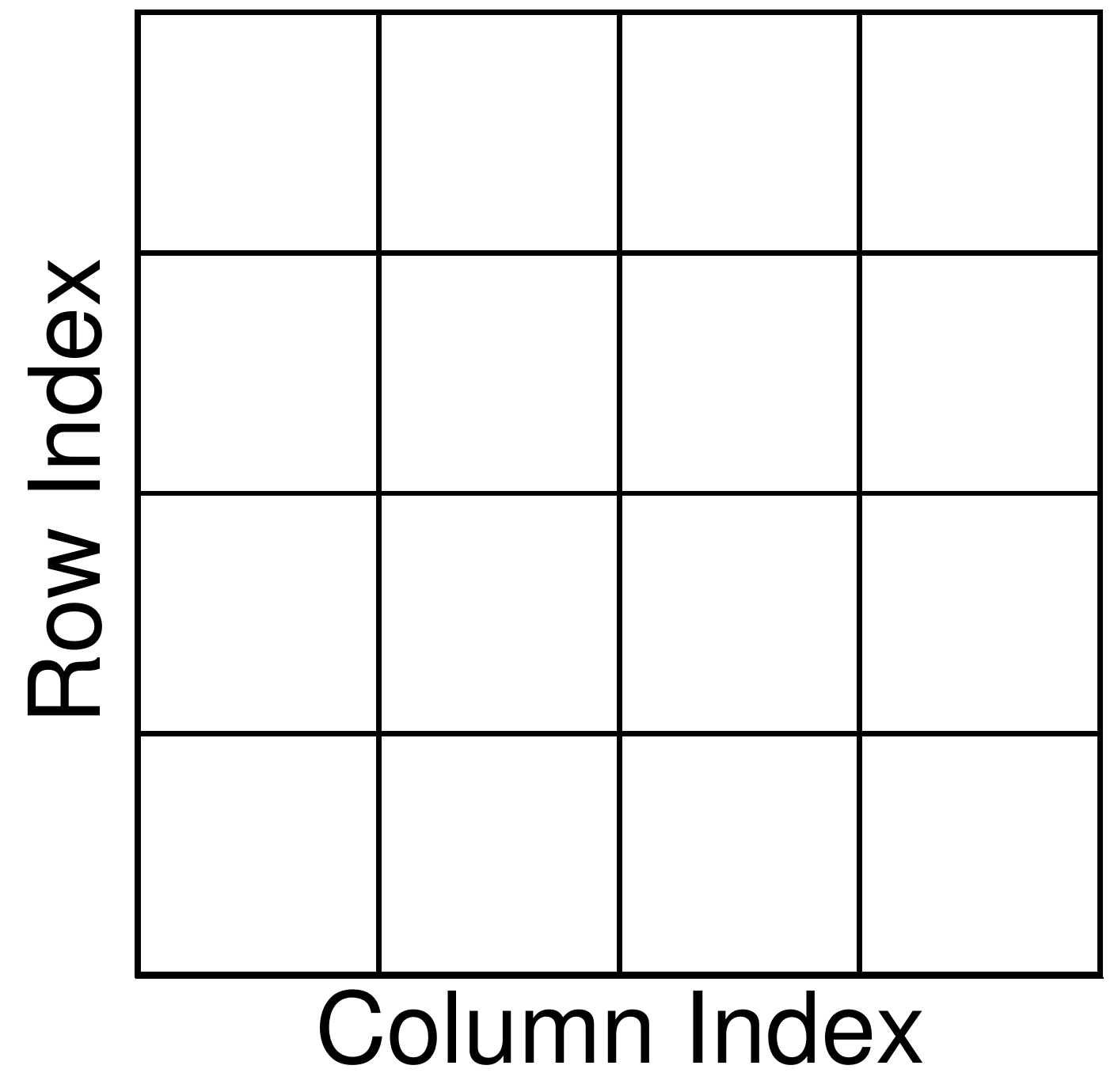}&
      \includegraphics[height=1.1in]{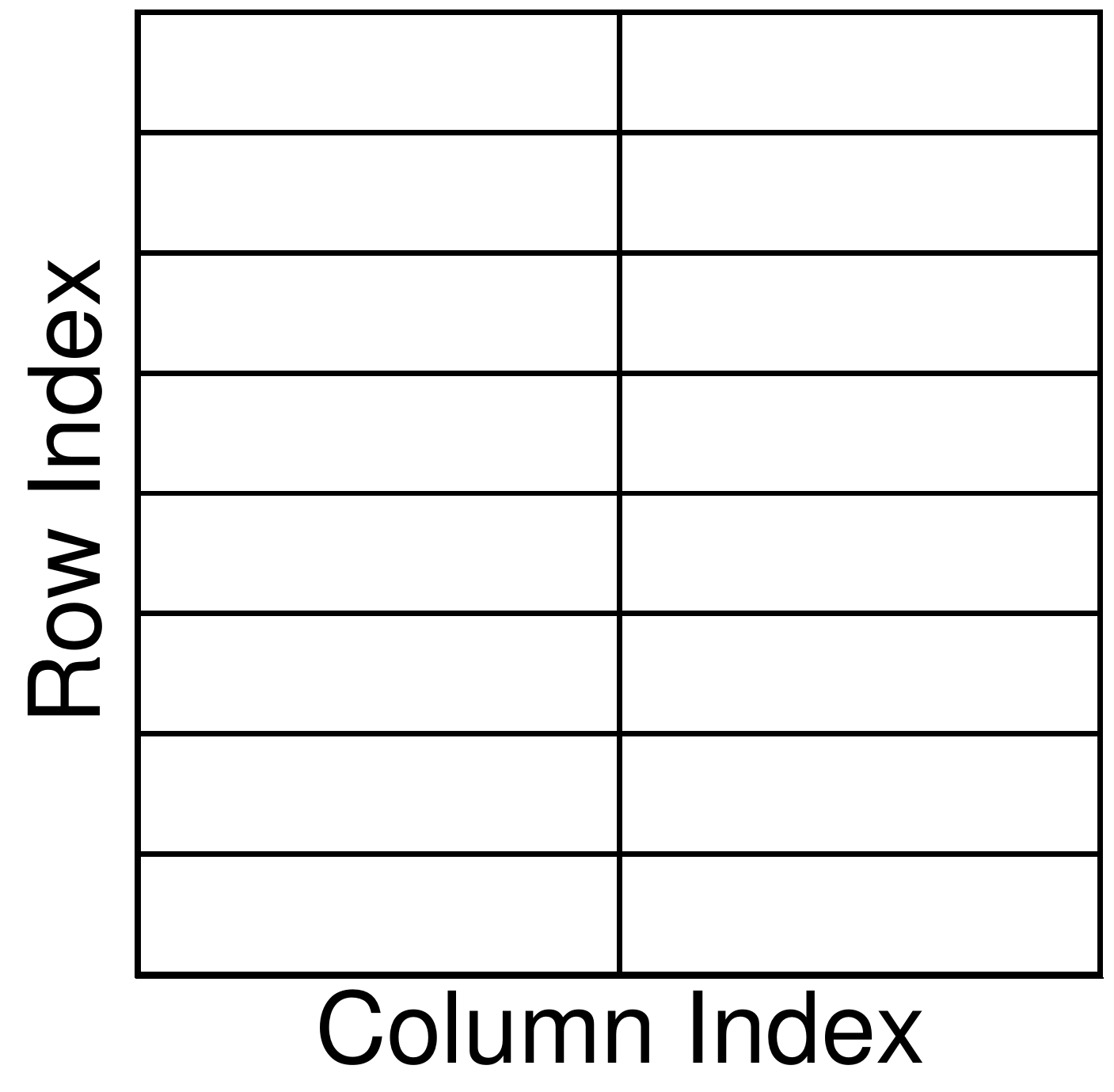}&
      \includegraphics[height=1.1in]{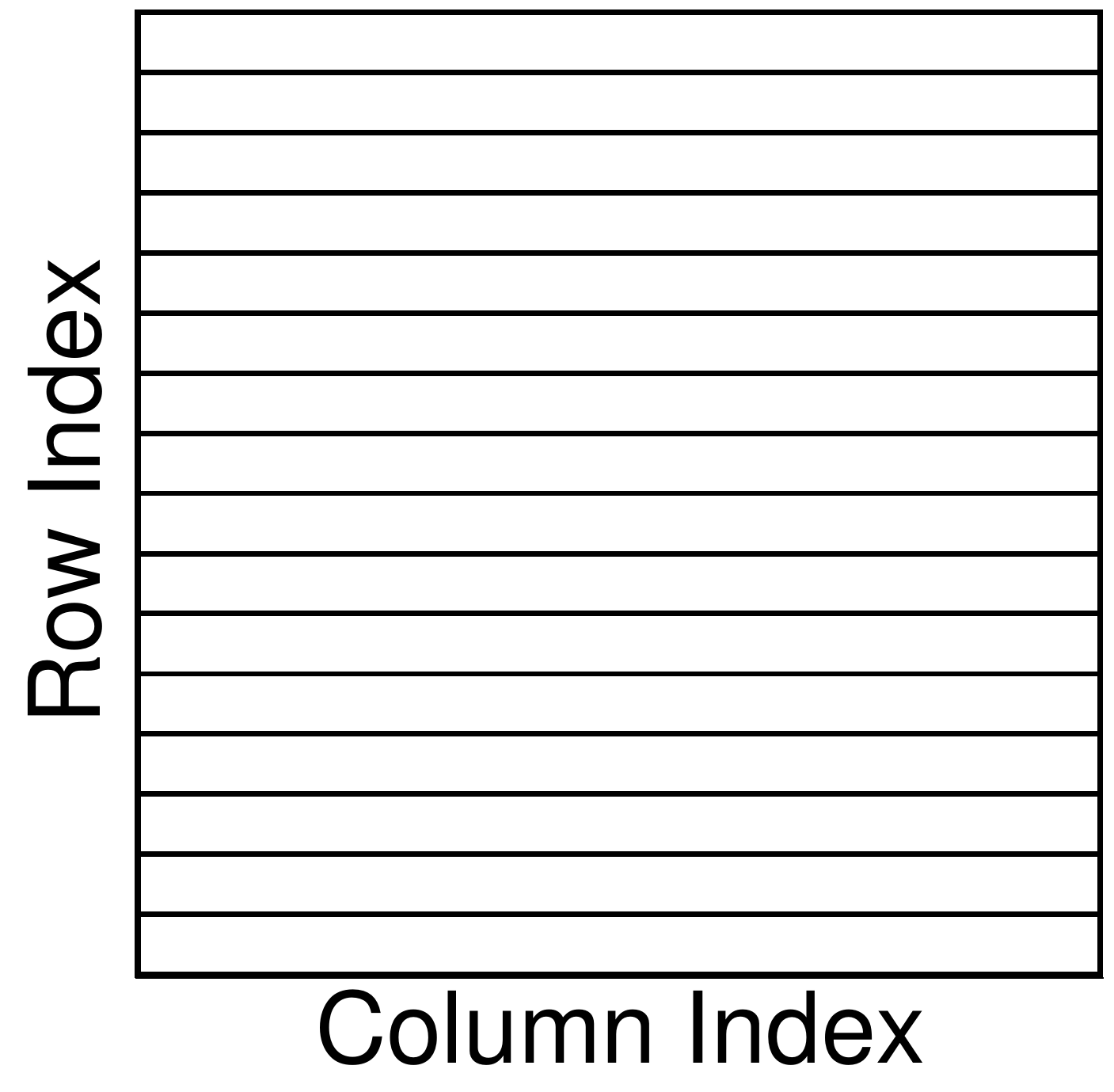}
    \end{tabular}
  \end{center}
  \caption{Hierarchical decomposition of the row and column indices
    of a one-dimensional complementary low-rank matrix of size
     $16\times 16$.  The trees $T_{X}$ ($T_{\Omega}$)
    has a root containing $16$ column (row) indices and leaves
    containing a single column (row) index.  The rectangles above
    indicate some of the low-rank submatrices.}
\label{fig:submatrices}
\end{figure}

$K:=e^{2\pi\i  UV^*}$ is a complementary low-rank matrix that has been widely studied in \cite{Butterfly3,Butterfly4,BF,MBF,Butterfly1,Butterfly2,
Butterfly5}. Let $X$ and $\Omega$ be the row and column index sets of $e^{2\pi\i  UV^*}$. Two trees $T_X$ and $T_\Omega$ of the same depth $L=O(\log N)$,
associated with $X$ and $\Omega$ respectively,
are constructed by dyadic partitioning. Denote the root level of the tree as level $0$ and
the leaf one as level $L$. Such a matrix $K$ of size $N\times N$
is said to satisfy the {\bf complementary low-rank property} if for
any level $\ell$, any node $A$ in $T_X$ at level $\ell$, and any node
$B$ in $T_\Omega$ at level $L-\ell$, the submatrix $K_{A,B}$, obtained
by restricting $K$ to the rows indexed by the points in $A$ and the
columns indexed by the points in $B$, is numerically low-rank.
See Figure \ref{fig:submatrices} for an illustration of low-rank submatrices
in a complementary low-rank matrix of size $16\times 16$.

In a special case when $K$ has an explicit formula, \cite{FIO09} proposed an $O(N\log N)$ butterfly algorithm to construct a data-sparse representation of $K$ using the low-rank factorizations of low-rank submatrices in the complementary low-rank structure. \cite{IBF} further optimized this algorithm and formulated it into the form of BF: 
\begin{equation}
\label{eqn:BFM}
K \approx  U^L G^{L-1}\cdots G^{h}  M^{h} (H^{h})^* \dots (H^{1})^* (V^0)^*,
\end{equation}
where the depth $L=\O(\log N)$ is assumed to be even,
$h=L/2$ is a middle level index, and all factors are sparse matrices with
$\O(N)$ nonzero entries. Storing and applying the BF requires only $O(N\log N)$ complexity. However, in a general case when only the low-rank factorization of the phase matrix $\Phi\approx UV^*$ is available, the state-of-the-art purely algebraic approach to construct the BF requires at least $O(N^{1.5})$ computational complexity \cite{BF}. Though the application of the BF is highly efficient, the precomputation of the factorization is still not practical when $N$ is large.

The IBF-MAT in this paper admits $O(N\log N)$ construction and application complexity, which would be a useful tool in developing nearly linear scaling algorithms to solve a wide class of differential and integral equations when incorporated into the schemes in \cite{LUBF,precon3,HSSBF,precon1,precon2}. The main difference between IBF-MAT and the BF in \cite{FIO09,IBF} is that, we apply Algorithm \ref{alg:ilr} in Section \ref{sec:LRF} to construct the low-rank factorization of low-rank submatrices, instead of the interpolative low-rank approximation in Section 2.1 in \cite{IBF}. Hence, to reduce the length of this paper, we only illustrate how Algorithm \ref{alg:ilr} in this paper is applied to design an $O(N\log N)$ butterfly algorithm. The reader is referred to \cite{IBF} for the routines that construct the data-sparse representation in the form of \eqref{eqn:BFM} using the new butterfly algorithm.

With no loss of generality, we assume that $K=e^{2\pi\i  UV^*}$ coming from the discretization of $K(x,\xi)=e^{2\pi\i  \Phi(x,\xi)}$ with a uniform grid. Given an input vector $\{f(\xi), \xi \in \Omega\}$, the goal is to compute the {\em potential vector} $\{g(x), x\in X\}$ defined by
\[
g(x) = \sum_{\xi\in \Omega} K(x,\xi) f(\xi), \quad x \in X.
\]

The main data structure of the butterfly
algorithm is a pair of dyadic trees $T_X$ and $T_\Omega$. Recall that
for any pair of intervals $A\times B \in T_X\times T_\Omega$
obeying the condition $\ell_A + \ell_B = L$,
 the submatrix $\{K(x,\xi)\}_{x \in A, \xi \in B}$ is
approximately of a constant rank. An explicit method to construct
its low-rank approximation is given by Algorithm \ref{alg:ilr}.
More precisely, for any $\eps>0$,
there exists a constant $r_\eps$ independent of $N$ and two sets of
functions $\{\alpha^{AB}_t (x)\}_{1\le t \le r_\eps}$ and
$\{\beta^{AB}_t (\xi)\}_{1\le t \le r_\eps}$  given in
\eqref{eqn:ab1} or \eqref{eqn:ab2} such that
\begin{equation}
  \left| K(x,\xi) - \sum_{t=1}^{r_\eps} \alpha^{AB}_t(x)
  \beta^{AB}_t(\xi) \right| \le \eps, \quad \forall x\in A, \forall
  \xi\in B.
  \label{eqn:glr}
\end{equation}

For a given interval $B$ in $\Omega$, define $u^B(x)$ to be the {\em
  restricted potential} over the sources $\xi\in B$
\[
u^B(x) = \sum_{\xi\in B} K(x,\xi) g(\xi).
\]
The low-rank property gives a compact expansion for $\{u^B(x)\}_{x\in
  A}$. Summing \eqref{eqn:glr} over $\xi\in B$ with coefficients $g(\xi)$
gives
\[
\left| u^B(x) - \sum_{t=1}^{r_\eps} \alpha^{AB}_t(x) \left( \sum_{\xi\in B} \beta^{AB}_t(\xi) g(\xi) \right) \right|
\le \left( \sum_{\xi\in B} |g(\xi)| \right) \eps,
\quad \forall x \in A.
\]
Therefore, if one can find coefficients $\{\coef{AB}{t}\}_{1\le t\le
  r_\eps}$ obeying
\begin{equation}
  \coef{AB}{t} \approx \sum_{\xi\in B} \beta^{AB}_t(\xi) g(\xi), \quad 1\le t\le r_\eps,
  \label{eqn:delta}
\end{equation}
then the restricted potential $\{u^B(x)\}_{x\in A}$ admits a compact
expansion
\[
\left| u^B(x) - \sum_{t=1}^{r_\eps} \alpha^{AB}_t(x) \coef{AB}{t} \right| \le \left( \sum_{\xi\in B} |g(\xi)| \right) \eps,
\quad \forall x\in A.
\]
%A key point of the butterfly algorithm is that for each pair $(A,B)$,
%he number of terms in the expansion is independent of $N$.
The butterfly algorithm below provides an efficient way for
computing $\{\coef{AB}{t}\}_{1\le t\le r_\eps}$ recursively.  The
general structure of the algorithm consists of a top-down traversal of
$T_X$ and a bottom-up traversal of $T_\Omega$, carried out
simultaneously. A schematic illustration of the data flow
in this algorithm is provided in Figure
\ref{fig:domain-tree-BA}.

\vspace{0.5cm}
\begin{algo}{Butterfly algorithm}
\label{aglo:BA}

\begin{enumerate}
\item{\emph{Preliminaries.}} Construct the trees $T_X$ and $T_\Omega$.

\item{\emph{Initialization.}} Let $A$ be the root of $T_X$. For each leaf interval $B$ of
  $T_\Omega$, construct the expansion coefficients $\{
  \coef{AB}{t}\}_{1\le t \le r_\eps}$ for the potential
  $\{u^B(x)\}_{x\in A}$ by simply setting
  \begin{equation}
    \coef{AB}{t} = \sum_{\xi\in B} \beta^{AB}_t(\xi) g(\xi), \quad 1\le t \le r_\eps.
  \label{eqn:bf1}
  \end{equation}
    By the interpolative low-rank approximation in Algorithm \ref{alg:ilr} applied to $e^{2\pi\i\Phi(A,B)}$ in the variable $\xi$ in $B$, we can define the expansion coefficients $\{\coef{AB}{t}\}_{1\le
    t\le r_\eps}$ by
  \begin{equation}
    \label{eqn:1}
    \coef{AB}{t} :=
    e^{-2\pi\i \Phi(c_A,\xi_t^B)} \sum_{\xi\in B}\left( M_t^B(\xi)
    e^{2\pi\i\Phi(c_A,\xi)}g(\xi)\right),
  \end{equation}
  where $\{\xi_t^B\}_{1\leq t\leq r_\epsilon}$ is the set of grid points adapted to $B$ by \eqref{eqn:gdB}.

  \item{\emph{Recursion.}} For $\ell = 1, 2, \ldots, L/2$, visit level $\ell$ in $T_X$ and
  level $L-\ell$ in $T_\Omega$. For each pair $(A,B)$ with $\ell_A =
  \ell$ and $\ell_B = L-\ell$, construct the expansion coefficients
  $\{\coef{AB}{t}\}_{1\le t \le r_\eps}$ for the potential
  $\{u^B(x)\}_{x\in A}$ using the low-rank representation constructed
  at the previous level. Let
  $P$ be $A$'s parent and $C$ be a child of $B$. Throughout, we shall
  use the notation $C\succ B$ when $C$ is a child of $B$. At level
  $\ell-1$, the expansion coefficients $\{\coef{PC}{s}\}_{1\le
    s\le r_\eps}$ of $\{u^{C}(x)\}_{x\in P}$ are readily available
  and we have
  \[
  \left| u^{C}(x) - \sum_{s=1}^{r_\eps} \alpha^{PC}_{s}(x) \coef{PC}{s} \right| \le \left( \sum_{\xi\in C} |g(\xi)| \right) \eps,
  \quad \forall x\in P.
  \]
  Since $u^B(x) = \sum_{C\succ B} u^{C}(x)$, the previous inequality
  implies that
  \[
  \left| u^B(x) - \sum_{C\succ B} \sum_{s=1}^{r_\eps} \alpha^{PC}_{s}(x) \coef{PC}{s} \right| \le \left( \sum_{\xi\in B} |g(\xi)| \right) \eps,
  \quad \forall x\in P.
  \]
  Since $A \subset P$, the above approximation is of course true for
  any $x \in A$. However, since $\ell_A + \ell_B = L$, the sequence of
  restricted potentials $\{u^B(x)\}_{x\in A}$ also has a low-rank
  approximation of size $r_\eps$, namely,
  \[
  \left| u^B(x) - \sum_{t=1}^{r_\eps} \alpha^{AB}_t(x) \coef{AB}{t} \right| \le \left( \sum_{\xi\in B} |g(\xi)| \right) \eps,
  \quad \forall x\in A.
  \]
  Combining the last two approximations, we obtain that
  $\{\coef{AB}{t}\}_{1\le t\le r_\eps}$ should obey
  \begin{equation}
    \sum_{t=1}^{r_\eps} \alpha^{AB}_t(x) \coef{AB}{t} \approx
    \sum_{C\succ B} \sum_{s=1}^{r_\eps} \alpha^{PC}_{s}(x) \coef{PC}{s}, \quad \forall x\in A.
    \label{eqn:bf2}
  \end{equation}
  This is an over-determined linear system for
  $\{\coef{AB}{t}\}_{1\le t\le r_\eps}$ when
  $\{\coef{PC}{s}\}_{1\le s\le r_\eps,C\succ B}$ are available.
  The butterfly algorithm uses an
  efficient linear transformation approximately mapping
  $\{\coef{PC}{s}\}_{1\le s\le r_\eps,C\succ B}$ into
  $\{\coef{AB}{t}\}_{1\le t\le r_\eps}$ as follows
  \begin{equation}
    \label{eqn:2}
    \coef{AB}{t} := e^{-2\pi\i\Phi(c_A,\xi_t^B)}\sum_{C\succ B}\sum_{s=1}^{r_\eps} M_t^B(\xi_{s}^{C})
    e^{2\pi\i\Phi(c_A,\xi_{s}^{C})}\coef{PC}{s},
  \end{equation}
  where $\{\xi_t^B\}_{1\leq t\leq r_\epsilon}$ (and $\{\xi_t^C\}_{1\leq t\leq r_\epsilon}$) is the set of grid points adapted to $B$ (and $C$) by \eqref{eqn:gdB}.
  
\item{\emph{Switch.}} For the levels visited, interpolation is applied in variable $\xi$, while the interpolation
  is applied in variable $x$ for levels $\ell>L/2$.  Hence, we
  are switching the interpolation variable in Algorithm \ref{alg:idrLR} at this step.  Now we are
  still working on level $\ell=L/2$ and the same domain pairs
  $(A,B)$ in the last step.  Let $\coef{AB}{s}$ denote the expansion
  coefficients obtained by interpolative low-rank factorization using Algorithm \ref{alg:ilr} applied to $e^{2\pi\i\Phi(A,B)}$ in the variable $\xi$ in $B$ in the last step.  Correspondingly, $\{\xi_s^B\}_s$ are the interpolation grid points in
  $B$ in the last step.  We take advantage of the interpolation in
  variable $x$ in $A$ using Algorithm \ref{alg:ilr} applied to $e^{2\pi\i\Phi(A,B)}$ and generate grid points $\{x_t^A\}_{1\le t\le
    r_\eps}$ in $A$ by \eqref{eqn:gdA}.  Then we can define new expansion coefficients
  \[
  \coef{AB}{t} := \sum_{s=1}^{r_\eps} e^{2\pi\i\Phi(x_t^A,\xi_s^B)}\coef{AB}{s}.
  \]
\item{\emph{Recursion.}} Similar to the discussion in Step $3$,
  we go up in tree $T_{\Omega}$ and down in tree
  $T_X$ at the same time until we reach the level
  $\ell=L$. We construct the low-rank approximation functions by interpolation in variable $x$ using Algorithm \ref{alg:ilr} as follows:
  \begin{eqnarray}
    \label{eqn:intx}
    \alpha_t^{AB}(x)= e^{2\pi\i\Phi(x,c_B)}M_t^A(x) e^{-2\pi\i\Phi(x_t^A,c_B)},
    &\beta_t^{AB}(\xi)=e^{2\pi\i\Phi(x_t^A,\xi)},
  \end{eqnarray}
    where $\{x_t^A\}_{1\leq t\leq r_\epsilon}$ is the set of grid points adapted to $A$ by \eqref{eqn:gdA}.
  
  Hence, the new expansion coefficients $\{\coef{AB}{t}\}_{1\le
    t\le r_\eps}$ can be defined as
  \begin{equation}
    \label{eqn:3}
    \coef{AB}{t} := \sum_{C\succ B} e^{2\pi\i\Phi(x_t^A,c_{C})}
    \sum_{s=1}^{r_\eps}\left(M_{s}^{P}(x_{t}^{A})
    e^{-2\pi\i\Phi(x_{s}^{P},c_{C})}\coef{PC}{s}\right),
  \end{equation}
  where again $P$ is $A$'s parent and $C$ is a child interval of $B$.
\item{\emph{Termination.}} Finally, $\ell = L$ and set $B$ to be the root node of
  $T_\Omega$. For each leaf interval $A \in T_X$, use the constructed
  expansion coefficients $\{\coef{AB}{t}\}_{1\le t\le r_\eps}$  in \eqref{eqn:3} to
  evaluate $u^B(x)$ for each $x \in A$,
  \begin{equation}
    \begin{split}
    u(x) = u^B(x)&= \sum_{t=1}^{r_\eps} \alpha^{AB}_t (x) \coef{AB}{t}\\
    &=e^{2\pi\i\Phi(x,c_B)} \sum_{t=1}^{r_\eps}\left(M_t^A(x)e^{-2\pi\i\Phi(x_t^A,c_B)} \coef{AB}{t} \right),
    \label{eqn:bf3}
    \end{split}
  \end{equation}
    where $\{x_t^A\}_{1\leq t\leq r_\epsilon}$ is the set of grid points adapted to $A$ by \eqref{eqn:gdA}.

\end{enumerate}

\end{algo}

\input{figure/fig-domain-tree-BA}

Like the butterfly algorithm in \cite{FIO09}, Algorithm $4.1$ is a direct approach that use the low-rank matrix factorization by Algorithm \ref{alg:ilr} on-the-fly to evaluate the oscillatory integral transform
\[
g(x)=\int e^{2\pi\i \Phi(x,\xi)} f(\xi)d\xi
\]
in $O(N\log N)$ operations without precomputation. If repeated applications of the integral transform to multiple functions $f$'s are required, it is more efficient to keep the low-rank matrix factorizations and arrange them into the form of BF in \eqref{eqn:BFM}. Besides, the rank provided by interpolative factorization is far from optimal, which motivates the structure-preseving sweeping matrix compression technique in \cite{IBF} to further compress the preliminary BF by interpolative factorization to obtain a sparser BF, which is the IBF-MAT of the kernel $e^{2\pi\i \Phi(x,\xi)}$ in this paper. The reader is referred to \cite{IBF} for a complete re-compression algorithm. 

%======================================================
\section{Numerical results}
\label{sec:results}

This section presents several numerical examples to demonstrate the efficiency of the proposed unified framework. 
The numerical results were obtained
on a computer with
Intel\textsuperscript{\textregistered}
Xeon\textsuperscript{\textregistered}
CPU X5690 @ 3.47GHz (6 core/socket)
and 128 GB RAM.
All implementations are in MATLAB\textsuperscript{\textregistered}
and are available per request. This new framework will be incorporated into the ButterflyLab\footnote{Available on \url{https://github.com/ButterflyLab}.} in the future.

Let $\{g^d(x),x\in X\}$, $\{g^b(x),x\in X\}$
and $\{g^n(x),x\in X\}$ denote the results
given by the direct matrix-vector multiplication,
IBF-MAT,
and NUFFT, respectively.
The accuracy of applying fast algorithms is estimated by the relative error defined as follows,
\begin{equation}
\epsilon^b = \sqrt{\cfrac{\sum_{x\in S}|g^b(x)-g^d(x)|^2}
{\sum_{x\in S}|g^d(x)|^2}}
\quad \text{and} \quad
\epsilon^n = \sqrt{\cfrac{\sum_{x\in S}|g^n(x)-g^d(x)|^2}
{\sum_{x\in S}|g^d(x)|^2}},
\end{equation}
where $S$ is an index set containing $256$ randomly sampled row indices of the kernel matrix $K$. The error for recovering the amplitude function is defined as
\begin{equation}
\epsilon^{amp} = \frac{\|\mathcal{A}(S,S)-U(S,:)V(:,S)^*\|_2}{\|\mathcal{A}(S,S)\|_2},
\end{equation}
where $\mathcal{A}$ is the amplitude matrix and $UV^*$ is its low-rank recovery. The error for recovering the phase and the kernel functions are defined similarly and denoted as $\epsilon^{pha}$ and $\epsilon^{K}$, respectively.

\subsection{Accuracy and scaling of low-rank matrix recovery and IBF-MAT}

In first part of the numerical section, we present numerical results of several examples to demonstrate the accuracy and asymptotic scaling of the proposed low-rank matrix recovery for amplitude and phase functions, and IBF-MAT. With no loss of generality, we only focus on Scenarios $1$ and $2$ of indirect access. For the first scenario, we apply the proposed algorithms to evaluate a Fourier integral operator (FIO) in 1D and a Hankel matrix transform. For the second scenario, we compute the IBF-MAT of the composition of two FIO's when we only have the BF representing each FIO. 

\vspace{0.5cm}
{\bf One-dimensional FIO}

Our first example is to evaluate a one-dimensional FIO \cite{BF} of the
following form:
\begin{equation}
\label{eq:example1}
g(x) = \int_{\mathbb{R}}\alpha(x,\xi)e^{2\pi \imath \Phi(x,\xi)}\widehat{f}(\xi)d\xi,
\end{equation}
where $\widehat{f}$ is the {Fourier} transform of $f$, $\alpha(x,\xi)=1$, 
and $\Phi(x,\xi)$ is a phase function given by
\begin{equation}
\Phi(x,\xi) = x\cdot \xi + c(x)|\xi|,~~~c(x) = (2+0.2\sin(2\pi x))/16.
\label{eqn:1D-FIO-kernal}
\end{equation}
The discretization of \eqref{eq:example1} is
\begin{equation}
\label{eqn:1D-FIO}
g(x_i) = \sum_{\xi_j}\alpha(x_i,\xi_j) e^{2\pi \imath \Phi(x_i,\xi_j)}\widehat{f}(\xi_j),
\quad i,j=1,2,\dots,N,
\end{equation}
where $\{x_i\}$ and $\{\xi_j\}$ are points uniformly distributed
in $[0,1)$ and $[-N/2,N/2)$ following
\begin{equation}\label{eqn:1D-xandxi}
x_i = (i-1)/N \text{ and } \xi_j = j-1-N/2.
\end{equation}

This example is for Scenario $1$ in Table \ref{tab:sc}. The unified framework is applied to recover the amplitude and phase functions in a form of low-rank matrix factorization, compute the IBF-MAT of the kernel function, and apply the IBF-MAT as in \eqref{eqn:tg} to a randomly generated $f$ in \eqref{eq:example1} to obtain $g$. Table~\ref{tab:1D-FIO} summarizes the results of this example for different grid sizes $N$ and numbers of interpolation points $r_\epsilon$. In the low-rank approximations of amplitude and phase functions, the rank parameter is $20$ and the over-sampling parameter is $5$.

\begin{table}[htp]
\centering
\begin{tabular}{rccccccccc}
\toprule
   $N,r_\epsilon$ & $\epsilon^b$ & $\epsilon^{K}$ &$\epsilon^{pha}$ &$\epsilon^{amp}$ & $T_{rec}(min)$& $T_{fac}(min)$
                             & $T_{app}(sec)$ & $T_d/T_{app}$\\
\toprule 
   1024, 6 & 2.52e-04 & 7.70e-11 & 7.70e-11 & 1.22e-15 & 1.80e-02 & 4.45e-02 &  6.17e-03 & 1.51e+01 \\
   1024, 8 & 2.60e-06 & 2.82e-12 & 2.82e-12 & 1.23e-15 & 1.07e-02 & 4.24e-02 &  3.91e-03 & 2.03e+01 \\
   1024,10 & 1.69e-08 & 3.16e-12 & 3.16e-12 & 1.22e-15 & 1.04e-02 & 3.35e-02 &  5.31e-03 & 1.72e+01 \\
   1024,12 & 6.21e-11 & 3.12e-12 & 3.12e-12 & 1.22e-15 & 1.08e-02 & 3.36e-02 &  3.95e-03 & 1.92e+01 \\
\toprule 
   4096, 6 & 3.38e-04 & 2.17e-11 & 2.17e-11 & 1.20e-15 & 4.06e-02 & 2.05e-01 &  1.32e-02 & 7.85e+01 \\
   4096, 8 & 3.16e-06 & 3.15e-11 & 3.15e-11 & 1.20e-15 & 4.21e-02 & 2.26e-01 &  1.62e-02 & 4.52e+01 \\
   4096,10 & 1.84e-08 & 6.67e-11 & 6.67e-11 & 1.31e-15 & 4.07e-02 & 1.83e-01 &  1.66e-02 & 4.38e+01 \\
   4096,12 & 7.87e-11 & 2.23e-11 & 2.23e-11 & 1.31e-15 & 4.06e-02 & 2.11e-01 &  2.34e-02 & 3.53e+01 \\
\toprule 
  16384, 6 & 3.87e-04 & 3.98e-10 & 3.98e-10 & 1.23e-15 & 1.53e-01 & 1.04e+00 &  4.78e-02 & 1.69e+02 \\
  16384, 8 & 3.98e-06 & 4.77e-11 & 4.77e-11 & 1.22e-15 & 1.54e-01 & 1.13e+00 &  7.39e-02 & 1.08e+02 \\
  16384,10 & 2.18e-08 & 2.64e-10 & 2.64e-10 & 1.22e-15 & 1.47e-01 & 9.73e-01 &  8.39e-02 & 1.02e+02 \\
  16384,12 & 1.87e-10 & 2.12e-10 & 2.12e-10 & 1.23e-15 & 1.47e-01 & 1.13e+00 &  1.14e-01 & 6.67e+01 \\
\toprule 
  65536, 6 & 4.85e-04 & 2.83e-09 & 2.83e-09 & 1.22e-15 & 5.81e-01 & 4.80e+00 &  1.96e-01 & 5.36e+02 \\
  65536, 8 & 5.35e-06 & 2.30e-09 & 2.30e-09 & 1.22e-15 & 5.77e-01 & 5.41e+00 &  3.07e-01 & 3.32e+02 \\
  65536,10 & 3.18e-08 & 3.77e-09 & 3.77e-09 & 1.22e-15 & 6.01e-01 & 5.07e+00 &  3.94e-01 & 2.91e+02 \\
  65536,12 & 2.01e-09 & 3.47e-09 & 3.47e-09 & 1.22e-15 & 5.96e-01 & 5.92e+00 &  5.40e-01 & 2.63e+02 \\
\toprule 
 262144, 6 & 5.55e-04 & 5.46e-09 & 5.46e-09 & 1.27e-15 & 2.32e+00 & 2.31e+01 &  8.80e-01 & 1.90e+03 \\
 262144, 8 & 4.51e-06 & 7.31e-09 & 7.31e-09 & 1.14e-15 & 2.32e+00 & 2.73e+01 &  1.48e+00 & 1.12e+03 \\
 262144,10 & 3.80e-08 & 2.23e-08 & 2.23e-08 & 1.25e-15 & 2.33e+00 & 2.51e+01 &  1.92e+00 & 8.69e+02 \\
 262144,12 & 7.70e-09 & 9.88e-09 & 9.88e-09 & 1.25e-15 & 2.33e+00 & 2.94e+01 &  2.55e+00 & 7.82e+02 \\
\bottomrule
\end{tabular}
\caption{Numerical results for the one-dimensional FIO given
  in \eqref{eqn:1D-FIO}. $T_{rec}$ is the time for recovering the amplitude and phase functions, $T_{fac}$ is the time for computing the IBF-MAT, $T_{app}$ is the time for applying the IBF-MAT, and $T_d$ is the time for a direct summation in \eqref{eqn:1D-FIO}.}
\label{tab:1D-FIO}
\end{table}

Table~\ref{tab:1D-FIO} shows that for a fixed number of interpolation points $r_\epsilon$, and a rank parameter for the amplitude and phase functions, the accuracy of the low-rank matrix recovery and the IBF-MAT stay in almost the same order, though the accuracy becomes slightly worse as the problem size increases. The slightly increasing error is due to the randomness of the proposed algorithm. As the problem size increases, the probability for capturing the low-rank matrix with a fixed rank parameter becomes smaller. Although the phase function is not smooth at $\xi=0$, the proposed algorithm is still able to recover the phase function accurately.

As for the computational complexity, both the factorization time and the application time of the IBF-MAT, and the reconstruction time of the amplitude and phase functions scales like $N\log N$. Every time we quadripule the problem size, the time increases on average by a factor of 4 to 6, and the increasing factor tends to decrease as the problem size increases. The speed-up factor over the direct method increases quickly and it is very significant when the problem size is large.

\vspace{0.5cm}
{\bf Special function transform}

Next, we provide an example of a special function transform. Following the standard notation, we denote the Hankel function of the first kind of order $m$ by $H_m^{(1)}$. When $m$ is an integer, $H_m^{(1)}$ has a singularity at the origin and a branch cut along the negative real axis. We are interested in evaluating the sum of Hankel functions over different orders,
\begin{equation}
\label{eqn:HK}
g(x_i)=\sum_{j=1}^NH^{(1)}_{j-1}(x_i)f_j,\quad i=1,2,\dots,N,
\end{equation}
which is analogous to expansion in orthogonal polynomials. The points $x_i$ are defined via the formula
\[
x_i=N+\frac{2\pi}{3}(i-1),
\]
which are bounded away from zero. The entries of the matrix in the above matvec can be explicitly calculated on-the-fly in $O(1)$ operations per entry using asymptotic formulas. The unified framework will work for many other orthogonal transforms in the oscillatory regime that admit smooth amplitude and phase functions. For more examples see \cite{James:2017,Bremer201815}.

This example is also for Scenario $1$ in Table \ref{tab:sc}. The unified framework is applied to recover the amplitude and phase functions in a form of low-rank matrix factorization, compute the IBF-MAT of the kernel function, and apply the IBF-MAT as in \eqref{eqn:tg} to a randomly generated $f$ to obtain $g$. Table~\ref{tab:HK} summarizes the results of this example for different grid sizes $N$ and numbers of interpolation points $r_\epsilon$. In the low-rank approximations of amplitude and phase functions, the rank parameter is $20$ and the over-sampling parameter is $5$.

\begin{table}[htp]
\centering
\begin{tabular}{rccccccccc}
\toprule
   $N,r_\epsilon$ & $\epsilon^b$ & $\epsilon^{K}$ &$\epsilon^{pha}$ &$\epsilon^{amp}$ & $T_{rec}(min)$& $T_{fac}(min)$
                             & $T_{app}(sec)$ & $T_d/T_{app}$\\
\toprule
   1024, 6 & 1.19e-04 & 3.07e-09 & 2.88e-09 & 6.78e-12 & 2.06e-02 & 3.65e-02 &  1.82e-02 & 3.79e+01 \\
   1024, 8 & 2.35e-06 & 5.77e-10 & 6.34e-10 & 1.88e-11 & 1.86e-02 & 3.61e-02 &  1.38e-02 & 6.29e+01 \\
   1024,10 & 2.26e-06 & 5.27e-10 & 5.75e-10 & 2.33e-11 & 1.84e-02 & 2.41e-02 &  1.23e-02 & 7.69e+01 \\
   1024,12 & 1.72e-07 & 5.21e-10 & 5.58e-10 & 1.70e-11 & 1.90e-02 & 2.61e-02 &  1.23e-02 & 7.27e+01 \\
\toprule 
   4096, 6 & 3.66e-05 & 1.73e-07 & 1.92e-07 & 2.85e-10 & 6.38e-02 & 1.80e-01 &  6.53e-02 & 1.41e+02 \\
   4096, 8 & 9.03e-06 & 2.52e-09 & 1.42e-09 & 1.16e-10 & 6.68e-02 & 1.98e-01 &  8.72e-02 & 9.98e+01 \\
   4096,10 & 1.97e-06 & 5.83e-09 & 3.16e-09 & 1.11e-10 & 6.66e-02 & 1.52e-01 &  8.35e-02 & 1.16e+02 \\
   4096,12 & 4.93e-07 & 9.66e-08 & 8.34e-08 & 3.05e-11 & 6.66e-02 & 1.64e-01 &  9.15e-02 & 1.01e+02 \\
\toprule 
  16384, 6 & 2.82e-03 & 5.00e-07 & 3.19e-07 & 7.10e-10 & 2.48e-01 & 9.51e-01 &  3.70e-01 & 2.92e+02 \\
  16384, 8 & 1.66e-04 & 7.16e-07 & 6.00e-07 & 9.17e-10 & 2.42e-01 & 1.02e+00 &  5.03e-01 & 2.04e+02 \\
  16384,10 & 4.21e-06 & 7.43e-08 & 3.75e-08 & 2.51e-09 & 2.49e-01 & 8.48e-01 &  4.66e-01 & 2.32e+02 \\
  16384,12 & 2.43e-07 & 3.61e-08 & 2.16e-08 & 1.87e-10 & 2.49e-01 & 8.88e-01 &  5.01e-01 & 2.08e+02 \\
\toprule 
  65536, 6 & 2.86e-03 & 2.51e-06 & 1.65e-07 & 3.97e-07 & 9.81e-01 & 4.56e+00 &  2.78e+00 & 6.65e+02 \\
  65536, 8 & 7.15e-06 & 2.98e-06 & 1.24e-06 & 1.11e-07 & 9.61e-01 & 4.96e+00 &  3.57e+00 & 4.74e+02 \\
  65536,10 & 8.50e-07 & 6.45e-06 & 3.40e-06 & 7.58e-11 & 9.59e-01 & 4.35e+00 &  4.04e+00 & 4.25e+02 \\
  65536,12 & 4.10e-05 & 1.99e-04 & 2.89e-06 & 3.12e-05 & 9.56e-01 & 4.67e+00 &  3.97e+00 & 4.69e+02 \\
\toprule 
 262144, 6 & 1.26e-03 & 3.82e-05 & 3.07e-05 & 9.52e-08 & 3.86e+00 & 2.22e+01 &  1.33e+01 & 1.81e+03 \\
 262144, 8 & 5.41e-06 & 9.77e-06 & 4.93e-06 & 1.26e-06 & 3.89e+00 & 2.42e+01 &  1.79e+01 & 1.38e+03 \\
 262144,10 & 1.01e-05 & 3.94e-05 & 1.61e-05 & 3.00e-06 & 3.90e+00 & 2.28e+01 &  2.37e+01 & 1.07e+03 \\
 262144,12 & 5.94e-05 & 1.58e-04 & 4.27e-06 & 2.17e-05 & 3.89e+00 & 2.44e+01 &  2.08e+01 & 1.17e+03 \\
\bottomrule
\end{tabular}
\caption{Numerical results for the Hankel function transform given
  in \eqref{eqn:HK}. $T_{rec}$ is the time for recovering the amplitude and phase functions, $T_{fac}$ is the time for computing the IBF-MAT, $T_{app}$ is the time for applying the IBF-MAT, and $T_d$ is the time for a direct summation in \eqref{eqn:HK}.}
\label{tab:HK}
\end{table}

\begin{figure}
  \begin{center}
    \begin{tabular}{c}
      \includegraphics[height=2in]{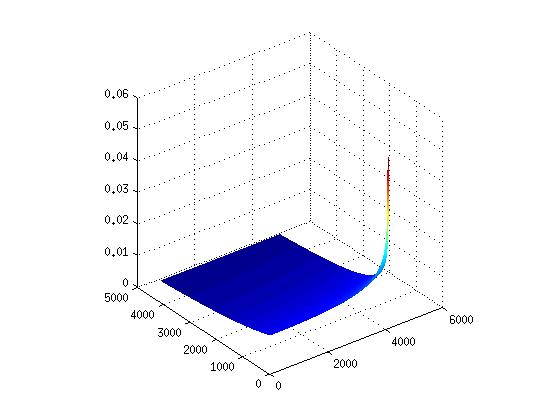}
    \end{tabular}
  \end{center}
  \caption{The exact amplitude function of the example in \eqref{eqn:HK}. There is a singularity point in the corner of the amplitude matrix, leading to an increasing rank of the amplitude matrix as the problem size increases.}
\label{fig:sing}
\end{figure}

The results in Table \ref{tab:HK} agree with the $\O(N\log N)$ complexity analysis and the speed-up factor over a direct summation is significant. The accuracy of the IBF-MAT becomes better if $r_\epsilon$ is larger and is almost independent of the problem size. Note that the recovery accuracy of the amplitude and phase functions becomes worse as $N$ increases. This is due to the fact that there is a singularity point in the corner of the amplitude matrix (see Figure \ref{fig:sing}), leading to an increasing rank of the amplitude matrix as the problem size increases. Besides, the randomized sampling algorithm in Algorithm \ref{alg:rlr} is not good in the presence of singularity, unless we know this singularity a prior so that we sample more at the corner. Hence, when $N>16384$ the accuracy of the low-rank amplitude and phase recovery is not very good and this influences the accuracy of the IBF-MAT, since the accuracy of the IBF-MAT is bounded below by the recovery accuracy. It is easy to fix this issue. After reconstructing the amplitude and phase, we can check singularity and reconstruct these functions again with adjusted sampling strategies to improve the accuracy. This works well in practice and we don't show the numerical results to save the space of the paper.

\vspace{0.5cm}
{\bf Composition of two FIO's in 1D}

The third example is to evaluate a composition of two FIO's of the
following form:
\begin{equation}
\label{eq:example4}
g(x)=\mathcal{L}\circ\mathcal{L}(f),
\end{equation}
where $\mathcal{L}$ is an FIO of the form
\begin{equation}
\label{eq:example3}
g(x) = \int_{\mathbb{R}}e^{2\pi \imath \Phi(x,\xi)}\widehat{f}(\xi)d\xi,
\end{equation}
where $\Phi(x,\xi)$ is a phase function given by
\begin{equation}
\Phi(x,\xi) = x\cdot \xi + c(x)\xi,~~~c(x) = (2+0.2\sin(2\pi x))/16.
\label{eqn:FIOK2}
\end{equation}
The discretization of \eqref{eq:example3} is similar to \eqref{eqn:1D-FIO}.

This example is for Scenario $2$ in Table \ref{tab:sc}. The unified framework is applied to recover the amplitude and phase functions in a form of low-rank matrix factorization, compute the IBF-MAT of the kernel function, and apply the IBF-MAT as in \eqref{eqn:tg} to a randomly generated $f$ in \eqref{eq:example4} to obtain $g$. Table~\ref{tab:2FIO} summarizes the results of this example for different grid sizes $N$ and numbers of interpolation points $r_\epsilon$. In the low-rank approximations of amplitude and phase functions, the rank parameter is $20$ and the over-sampling parameter is $5$.

We would like to emphasize that the composition of two FIO's results in an FIO with a phase function that is very singular at the point $\xi=0$. This leads to large-rank submatrices in the kernel matrix. In this case, we can adopt the multiscale butterfly algorithm/factorization in \cite{FIO14,MBF} to deal with this singularity. We have implemented the multiscale version of the IBF-MAT and present its numercial performance in Table \ref{tab:2FIO}. For the purpose of reducing the length of this paper, we don't introduce the multiscale IBF-MAT. The reader is referred to \cite{FIO14,MBF} for detailed description of the multiscale idea.

\begin{table}[h]
\centering
\begin{tabular}{rccccccccc}
\toprule
   $N,r_\epsilon$ & $\epsilon^b$ & $T_{rec}(min)$ & $T_{fac}(min)$
                             & $T_{app}(sec)$ \\
\toprule 
   1024, 6 & 3.13e-04 & 6.65e-02 & 2.70e-02 & 3.41e-02  \\
   1024, 8 & 3.65e-06 & 4.52e-02 & 2.75e-02 & 2.99e-02  \\
   1024,10 & 3.07e-08 & 4.55e-02 & 1.81e-02 & 3.65e-02  \\
   1024,12 & 4.25e-10 & 4.49e-02 & 2.04e-02 & 3.51e-02  \\
\toprule 
   4096, 6 & 3.94e-04 & 2.43e-01 & 1.60e-01 & 1.25e-01  \\
   4096, 8 & 4.59e-06 & 2.41e-01 & 1.81e-01 & 1.79e-01  \\
   4096,10 & 3.48e-08 & 2.47e-01 & 1.49e-01 & 2.54e-01  \\
   4096,12 & 9.24e-10 & 2.45e-01 & 1.71e-01 & 3.14e-01  \\
\toprule 
  16384, 6 & 4.58e-04 & 1.51e+00 & 8.92e-01 & 7.02e-01  \\
  16384, 8 & 5.42e-06 & 1.80e+00 & 1.02e+00 & 1.69e+00  \\
  16384,10 & 3.84e-08 & 1.72e+00 & 9.42e-01 & 1.70e+00  \\
  16384,12 & 1.69e-09 & 1.80e+00 & 1.08e+00 & 1.86e+00  \\
\toprule 
  65536, 6 & 5.22e-04 & 9.33e+00 & 4.61e+00 & 7.90e+00  \\
  65536, 8 & 6.29e-06 & 1.01e+01 & 5.36e+00 & 1.42e+01  \\
  65536,10 & 4.56e-08 & 9.54e+00 & 5.12e+00 & 2.47e+01  \\
  65536,12 & 9.25e-09 & 1.01e+01 & 5.72e+00 & 2.68e+01  \\
\toprule 
 262144, 6 & 5.86e-04 & 3.32e+01 & 2.22e+01 & 5.30e+01  \\
 262144, 8 & 7.16e-06 & 3.29e+01 & 2.51e+01 & 8.35e+01  \\
 262144,10 & 6.07e-08 & 3.28e+01 & 2.55e+01 & 1.40e+02  \\
 262144,12 & 2.45e-08 & 3.25e+01 & 3.00e+01 & 1.43e+02  \\
\bottomrule
\end{tabular}
\caption{Numerical results for the composition of two FIO's given
  in \eqref{eq:example4}. $T_{rec}$ is the time for recovering the amplitude and phase functions, $T_{fac}$ is the time for computing the multiscale IBF-MAT, and $T_{app}$ is the time for applying the multiscale IBF-MAT.}
\label{tab:2FIO}
\end{table}

Table~\ref{tab:2FIO} shows that for a fixed number of interpolation points $r_\epsilon$, and a rank parameter for the amplitude and phase functions, the accuracy of the low-rank matrix recovery and the multiscale IBF-MAT stay in almost the same order, though the accuracy becomes slightly worse as the problem size increases. The slightly increasing error is due to the randomness of the proposed algorithm as explained previously. There is no explicit formula for the amplitude and phase functions in this example. Hence, we cannot estimate the accuracy of the recovery algorithm. Since the accuracy of the multiscale IBF-MAT is bounded below by the accuracy of amplitude and phase recovery. We see that the recovery accuracy should be very good.

As for the computational complexity, both the factorization time and the application time of the IBF-MAT, and the reconstruction time of the amplitude and phase functions scales like $N\log N$. On average, when we quadripule the problem size, the time increases on average by a factor of 4 to 6, and the increasing factor tends to decrease as the problem size increases.

\subsection{Comparison of NUFFT and BF}

In the second part of the numerical section, we illustrate the $O(N)$ strategy in Algorithm \ref{alg:dc} for deciding whether we can use NUFFT in the oscillatory integral transform. {We will show that once the NUFFT is applicable, it is more efficient than the BF considering that the prefactor of the factorization and application time of the BF is larger than that of the NUFFT approach, when we require an approximate matvec with a high accuracy, no matter how many vectors in the matvec.} To this end, we will provide an example of FIO's in solving wave equations. {In the case of low accuracy requirement, according to the comparison of BF and NUFFT in Table $1$ and $2$ in \cite{IBF}, our conclusion just above still valid.}

Fast algorithms for solving wave equations with variable coefficients based on FIO's have been studied based on either the BF  in \cite{Yingwave} or the wave packet representation of the FIO's  in \cite{FIOwp2,wpFIO}. \cite{Yingwave} also proposed an approach to solve wave equations based on a carefully desgined NUFFT according to the explicit formulas of FIO's inspired by the work in \cite{FIO07}.  

We propose to apply the new NUFFT approach with dimension lifting for the evaluation of FIO's in solving wave equations. This new method does not rely on the explicit formula of an FIO and can be applied to more general scenarios. Besides, the dimension lifting idea could lead to fewer applications of the NUFFT, since the rank $r_\epsilon$ in \eqref{eqn:tg3} could be smaller compared to the NUFFT approach in \cite{Yingwave}. We will only provide a one-dimensional wave equation as an example to compare the performance of the new NUFFT approach and the BF approach for the evaluation of FIO's in solving wave equations. The application of the new NUFFT approach to solve higher dimensional wave equations will be left as a future work. 

In more particular, we solve the one-dimensional wave equation as follows:
\begin{eqnarray}\label{eqn:weqn}
\begin{cases}
\partial_{tt}u(x,t)-\partial_x(c^2(x)\partial_x u(x,t))=0 &\quad t>0,x\in[0,1)\\
u(x,0)=u_0(x)&\quad x\in[0,1)\\
\partial_t u(x,0)=u_1(x)& \quad x\in[0,1),
\end{cases}
\end{eqnarray}
where the boundary conditions are taken to be periodic. The theory of FIO's states that for a given smooth and positive $c(x)$ there exists a time $t^*$ that depends only on $c(x)$ such that for any $t<t^*$, the general solution of \eqref{eqn:weqn} is given by a summation of two FIO's:
\[
u(x,t)=\sum_{\xi\in\mathbb{Z}}e^{2\pi\i  \Phi_{\pm}(x,\xi,t)} \alpha_{\pm}(x,\xi,t)\widehat{f_{\pm}}(\xi),
\]
where $f_{\pm}$ are two functions depending on the initial conditions.

In this example, we assume that  $c(x)=2+\sin(2\pi x)$ and follow the ideas in \cite{Yingwave} to construct the FIO's in the solution operator of \eqref{eqn:weqn}. Without loss of generality, we focus on the evaluation of the  FIO
\begin{equation}\label{eqn:FIOwq}
\sum_{\xi\in\mathbb{Z}}e^{2\pi\i  \Phi_{+}(x,\xi,t)} \widehat{f_{+}}(\xi).
\end{equation}
The phase function $\Phi_+(x,\xi,t)$ satisfies the Hamiltonian-Jacobi equation
\begin{eqnarray}\label{eqn:HJ1}
\begin{cases}
\partial_t\Phi_+(x,\xi,t)-c(x)|\partial_x\Phi_+(x,\xi,t)|=0\\
\Phi_+(x,\xi,0)=x\cdot \xi.
\end{cases}
\end{eqnarray}
Note that $\Phi_+(x,\xi,t)$ is homogeneous of degree $1$ in $\xi$,  i.e., $\Phi_+(x,\lambda\xi,t)=\lambda \Phi_+(x,\xi,t)$ for $\lambda>0$. Therefore, we only need to evaluate $\Phi_+(x,\xi,t)$ at $\xi=\pm 1$. From the algebraic point of view, the phase matrix is piecewise rank-$1$, i.e., 
\begin{equation}\label{eqn:phlr}
\Phi_+(x,\xi,t)=
\begin{cases}
\Phi_+(x,1,t)\xi,\quad \forall \xi\geq 0,\\
-\Phi_+(x,-1,t)\xi,\quad\forall\xi<0.
\end{cases}
\end{equation}
In fact, to make the boundary condition periodic in $x$, $\Psi_+(x,\xi,t):=\Phi_+(x,\xi,t)-x\xi$ is introduced for $\xi=\pm 1$. Then we have
\begin{eqnarray}\label{eqn:HJ2}
\begin{cases}
\partial_t\Psi_+(x,\xi,t)-c(x)|\partial_x\Psi_+(x,\xi,t)+\xi |=0,\\
\Psi_+(x,\xi,0)=0.
\end{cases}
\end{eqnarray}
When $c(x)$ is a band-limited function, $\Psi_+(x,\xi,t)$ is a smooth function in $x$ when $t$ is sufficiently smaller than $t^*$. Hence, a small grid in $x$ is enough to discretize \eqref{eqn:HJ2}. The value of $\Psi_+$ on a finer grid in $x$ can evaluated by spectral interpolation using FFT. 

In the numerical examples here, we adopt a uniform grid with $512$ grid points for $x$ in $[0,1)$, and a time step size $\frac{1}{4096}$ to solve \eqref{eqn:HJ2}. The standard local Lax-Friedrichs Hamiltonian method is applied for $x$ and the third order TVD Runge-Kutta method is used for $t$ to solve \eqref{eqn:HJ2}. We vary the problem size $N$ of the evaluation in \eqref{eqn:FIOwq} and discretize $\Phi_+(x,\xi,t)$ with a uniform spacial grid with a step size $\frac{1}{N}$ for $x\in[0,1)$ and a uniform frequency grid with a step size $1$ for $\xi\in[-\frac{N}{2},\frac{N}{2})$.

By \eqref{eqn:phlr}, we solve \eqref{eqn:HJ2} and obtain a low-rank factorization of the phase matrix and apply IBF-MAT to evaluate \eqref{eqn:FIOwq}. Note that the phase matrix is piecewise rank-$1$, we can split the summation in \eqref{eqn:FIOwq} into two parts:
\begin{equation}\label{eqn:FIOwq2}
\sum_{\xi\in\{-\frac{N}{2},\dots,-1\}}e^{2\pi\i  \Phi_{+}(x,\xi,t)} \widehat{f_{+}}(\xi)+\sum_{\xi\in\{0,1,\cdots,\frac{N}{2}-1\}}e^{2\pi\i  \Phi_{+}(x,\xi,t)} \widehat{f_{+}}(\xi),
\end{equation}
and apply the one-dimensional NUFFT approach to evaluate the two summations in \eqref{eqn:FIOwq2}. Or we can also apply the two-dimensional NUFFT approach to compute the summation in \eqref{eqn:FIOwq}. 
%\footnote{For higher dimensional wave equations, the piecewise rank-$1$ property of the phase matrix is approximately valid and the number of pieces depends on the number of grid points of $\xi$ on the unit sphere used to solve \eqref{eqn:HJ2}. Since the phase function is smooth except for $\xi=0$, a small number of grid points on the unit sphere is enough, and hence the phase matrix is piecewise rank-$1$ with a small number of pieces. Therefore, we only need $O(1)$ NUFFT's even for high dimensional wave equations. This will be left as future work. }
The numerical results are summarized in Table \ref{tab:weq} and Table \ref{tab:weq2}. To make the accuracy of the BF and the NUFFT approaches comparable, we choose the rank parameter $r_\epsilon$ in the IBF-MAT as $12$, the accuracy tolerance $\epsilon$ in the IBF-MAT and the NUFFT as $1e-12$. 

%To show the performance of the IBF-MAT when it has competitive application time over the NUFFT approach, we set $r_\epsilon=3$ and $\epsilon=1e-12$ in the experiments in Table \ref{tab:weq2}. 

\begin{table}[h]
\centering
\begin{tabular}{rccccccccc}
\toprule
   $N$ & $t$ & $T_{FFT}(sec)$ &$T^b_{fac}(sec)$ &$T^b_{app}(sec)$ &$\epsilon^b$ & $T_{fac}^n(sec)$& $T_{app}^n(sec)$
                             & $\epsilon^n$ \\
\toprule
   1024 & 2.441e-04 & 2.09e-04 & 2.22e+00 & 2.90e-03 & 9.86e-13& 3.96e-03& 2.55e-03& 1.69e-13 \\
   1024 & 1.953e-03 & 2.09e-04 & 1.72e+00 & 1.97e-03 & 9.43e-13& 1.80e-03& 1.12e-03& 9.85e-14 \\
   1024 & 1.562e-02 & 2.09e-04 & 1.71e+00 & 1.99e-03 & 1.37e-12& 2.28e-03& 1.03e-03& 8.25e-14 \\
\toprule 
   4096 & 2.441e-04 & 2.07e-04 & 1.08e+01 & 1.14e-02 & 1.33e-12& 9.26e-04& 3.44e-03& 4.22e-13 \\
   4096 & 1.953e-03 & 2.07e-04 & 1.02e+01 & 1.12e-02 & 1.25e-12& 8.13e-04& 3.44e-03& 2.74e-13 \\
   4096 & 1.562e-02 & 2.07e-04 & 1.03e+01 & 1.14e-02 & 1.52e-12& 7.78e-04& 3.36e-03& 1.93e-13 \\
\toprule 
  16384 & 2.441e-04 & 3.21e-04 & 5.56e+01 & 5.63e-02 & 6.46e-12& 8.95e-04& 1.26e-02& 1.02e-12 \\
  16384 & 1.953e-03 & 3.21e-04 & 5.53e+01 & 8.04e-02 & 6.13e-12& 9.30e-04& 1.36e-02& 1.48e-12 \\
  16384 & 1.562e-02 & 3.21e-04 & 5.61e+01 & 5.65e-02 & 5.29e-12& 1.02e-03& 1.48e-02& 1.06e-12 \\
\toprule 
  65536 & 2.441e-04 & 2.91e-03 & 2.92e+02 & 2.70e-01 & 3.06e-12& 9.89e-04& 5.33e-02& 7.89e-12 \\
  65536 & 1.953e-03 & 2.91e-03 & 2.93e+02 & 2.72e-01 & 3.96e-12& 9.62e-04& 5.39e-02& 5.55e-12 \\
  65536 & 1.562e-02 & 2.91e-03 & 2.93e+02 & 3.14e-01 & 3.78e-12& 1.12e-03& 6.05e-02& 4.30e-12 \\
\toprule 
 262144 & 2.441e-04 & 5.41e-03 & 1.46e+03 & 1.23e+00 & 3.61e-12& 8.45e-04& 2.04e-01& 6.33e-12 \\
 262144 & 1.953e-03 & 5.41e-03 & 1.48e+03 & 1.42e+00 & 9.87e-12& 1.08e-03& 2.16e-01& 4.12e-11 \\
 262144 & 1.562e-02 & 5.41e-03 & 1.56e+03 & 1.31e+00 & 1.04e-11& 1.20e-03& 2.00e-01& 4.04e-11 \\
\bottomrule
\end{tabular}
\caption{Numerical results for the evaluation of \eqref{eqn:FIOwq2} for different problem sizes $N$ at different time $t$. $T_{FFT}$ is the runtime of a FFT on a vector of length $N$ as comparison. $T^b_{fac}$, $T^b_{app}$, $T_{fac}^n$, and $T_{app}^n$ are the factorization time and the application time for the IBF-MAT and the one-dimensional NUFFT, respectively. $\epsilon^b$ and $\epsilon^n$ are the relative evaluation error by the IBF-MAT and the NUFFT approaches, respectively.}
\label{tab:weq}
\end{table}

Numerical results in Table \ref{tab:weq} show that both the IBF-MAT and the one-dimensional NUFFT approach without dimension lifting admit $O(N\log N)$ factorization and application time. For almost the same evaluation accuracy, the one-dimensional NUFFT approach has a much smaller prefactor (about $O(1000)$ times smaller considering the total cost) making it more preferable. 

Numerical results in Table \ref{tab:weq2} show that the two-dimensional NUFFT approach with dimension lifting also admits $O(N\log N)$ factorization and application time. Though the BF might be a few times more efficient in some cases in terms of the application time, the NUFFT approach is still more preferable considering the expensive factorization time of the BF. Although the two-dimensional NUFFT approach is more expensive than the one-dimensional NUFFT method, the two-dimensional NUFFT approach doesn't rely on the piecewise rank-$1$ property of the phase function, and therefore is applicable in more general situations. 

Although we know that the NUFFT approach is applicable for \eqref{eqn:FIOwq} and \eqref{eqn:FIOwq2}, we still apply Algorithm \ref{alg:dc} to test its time scaling. The results of $T_{DEC}$ in Table \ref{tab:weq2} also verify that Algorithm \ref{alg:dc} for deciding whether we can apply the NUFFT approach has a linear scaling.

%======================================================
\section{Conclusion and discussion}
\label{sec:conclusion}

This paper introduced a unified framework for $O(N\log N)$ evaluation of the oscillatory integral transform $g(x) = \int \alpha(x,\xi)e^{2\pi\i \Phi(x,\xi)} f(\xi)d\xi$. This framework works for two cases: 1) explicit formulas for the amplitude and phase functions are known; 2) only indirect access of the amplitude and phase functions are available. In the case of indirect access, this paper proposed a novel fast algorithms for recovering the amplitude and phase functions in $O(N\log N)$ operations. Second, a new algorithm for the oscillatory integral transform based on the NUFFT and a dimension lifting technique is proposed. Finally, a new BF, the IBF-MAT, for amplitude and phase matrices in a form of a low-rank factorization is proposed. These two algorithms both requires only $O(N\log N)$ operations to evaluate the oscillatory integral transform.

\begin{table}[h]
\centering
\begin{tabular}{rccccccccc}
\toprule
   $N$ & $t$ & $T_{dec}(sec)$ &$T^b_{fac}(sec)$ &$T^b_{app}(sec)$ &$\epsilon^b$ & $T_{fac}^n(sec)$& $T_{app}^n(sec)$
                             & $\epsilon^n$ \\
\toprule
   1024 & 2.441e-04 & 1.33e-02 & 2.22e+00 & 2.90e-03 & 9.86e-13& 1.65e-03& 7.05e-04& 1.45e-13 \\
   1024 & 1.953e-03 & 1.63e-02 & 1.72e+00 & 1.97e-03 & 9.43e-13& 2.90e-03& 4.94e-03& 1.03e-13 \\
   1024 & 1.562e-02 & 1.43e-02 & 1.71e+00 & 1.99e-03 & 1.37e-12& 4.32e-04& 4.98e-03& 9.55e-14 \\
\toprule 
   4096 & 2.441e-04 & 4.09e-02 & 1.08e+01 & 1.14e-02 & 1.33e-12& 7.83e-04& 1.96e-03& 4.66e-13 \\
   4096 & 1.953e-03 & 3.83e-02 & 1.02e+01 & 1.12e-02 & 1.25e-12& 9.56e-04& 1.53e-02& 6.70e-13 \\
   4096 & 1.562e-02 & 3.87e-02 & 1.03e+01 & 1.14e-02 & 1.52e-12& 1.66e-03& 2.21e-02& 8.89e-13 \\
\toprule 
  16384 & 2.441e-04 & 1.40e-01 & 5.56e+01 & 5.63e-02 & 6.46e-12& 6.30e-04& 7.41e-03& 1.07e-12 \\
  16384 & 1.953e-03 & 1.80e-01 & 5.53e+01 & 8.04e-02 & 6.13e-12& 4.69e-04& 7.84e-02& 5.01e-12 \\
  16384 & 1.562e-02 & 1.47e-01 & 5.61e+01 & 5.65e-02 & 5.29e-12& 4.16e-04& 1.54e-01& 1.62e-12 \\
\toprule 
  65536 & 2.441e-04 & 6.60e-01 & 2.92e+02 & 2.70e-01 & 3.06e-12& 7.04e-04& 3.25e-02& 4.69e-12 \\
  65536 & 1.953e-03 & 6.54e-01 & 2.93e+02 & 2.72e-01 & 3.96e-12& 4.75e-04& 3.75e-01& 1.93e-11 \\
  65536 & 1.562e-02 & 6.76e-01 & 2.93e+02 & 3.14e-01 & 3.78e-12& 5.01e-04& 9.34e-01& 3.92e-11 \\
\toprule 
 262144 & 2.441e-04 & 2.94e+00 & 1.46e+03 & 1.23e+00 & 3.61e-12& 2.34e-03& 1.18e-01& 2.34e-11 \\
 262144 & 1.953e-03 & 3.29e+00 & 1.48e+03 & 1.42e+00 & 9.87e-12& 7.74e-04& 2.46e+00& 2.62e-10 \\
 262144 & 1.562e-02 & 3.15e+00 & 1.56e+03 & 1.31e+00 & 1.04e-11& 4.51e-04& 8.13e+00& 3.02e-10 \\
\bottomrule
\end{tabular}
\caption{Numerical results for the evaluation of \eqref{eqn:FIOwq} for different problem sizes $N$ at different time $t$. $T_{dec}$ is the runtime of Algorithm \ref{alg:dc}.  $T^b_{fac}$, $T^b_{app}$, $T_{fac}^n$, and $T_{app}^n$ are the factorization time and the application time for the IBF-MAT and the two-dimensional NUFFT approach by dimension lifting, respectively. $\epsilon^b$ and $\epsilon^n$ are the relative evaluation error by the IBF-MAT and the NUFFT approaches, respectively.}
\label{tab:weq2}
\end{table}

This unified framework would be very useful in develping efficient tools for fast special function transforms, solving wave equations, and solving electromagnetic (EM)  scattering problems. We have provided several examples to support these applications. For example, the state-of-the-art fast algorithm for computing the compositions of FIO's, which could be applied as a preconditioner for certain classes of parabolic and hyperbolic equations \cite{precon3,precon1,precon2}; a fast algorithm for solving wave equation via FIO's. We have explored the potential applications of the proposed framework to: 1) fast evaluation of other special functions \cite{James:2017,Bremer201815} to develop nearly linear scaling polynomial transforms; 2) fast solvers developed in \cite{LUBF,HSSBF} for nearly linear algorithms for solving high-frequency EM equations. Numerical results will be reported in forthcoming papers.

{\bf Acknowledgments.}  The author thanks the fruitful discussion with Lexing Ying and the support of the start-up package at the National University of Singapore.

\bibliographystyle{abbrv}
\bibliography{ref}

\end{document}

%% file: figure/fig-domain-tree-BA.tex
\begin{figure}
\centering
\begin{tikzpicture}[scale=0.3]

\draw [->,double,double distance=1pt] (-6,9) -> (-6,1);

\draw (0,10) -- (-4,0) -- (4,0) -- cycle;

\draw [->,double,double distance=1pt] (20,1) -> (20,9);

\draw (14,10) -- (10,0) -- (18,0) -- cycle;

\draw [latex-latex,line width=1pt] (9,5) -> (5,5);
\draw [latex-latex,line width=1pt] (9,9) -> (5,1);
\draw [latex-latex,line width=1pt] (9,1) -> (5,9);

\coordinate [label=above:$T_X$] (X) at (0,10);
\coordinate [label=above:$T_\Omega$] (X) at (14,10);

\draw (-2,5) -- (2,5);
\draw (12,5) -- (16,5);

\coordinate [label=right:$\frac{L}{2}$] (X) at (16,5);
\coordinate [label=left:$\frac{L}{2}$] (X) at (-2,5);
\end{tikzpicture}
\caption{Trees of the row and column indices.
    Left: $T_X$ for the row indices $X$.
    Right: $T_\Omega$ for the column indices $\Omega$.
    The interaction between $A\in T_X$ and $B\in T_\Omega$
    starts at the root of $T_X$ and the leaves of $T_\Omega$. }
\label{fig:domain-tree-BA}
\end{figure}